\documentclass[12pt]{amsart}
\usepackage{graphicx}
\usepackage{amssymb}
\usepackage[all]{xy}
\numberwithin{equation}{section}

\newtheorem*{thmmain}{\bf Theorem}  
\newtheorem*{claim}{\bf Claim}      

\newtheorem{thm}[equation]{Theorem}

\newtheorem{lem}[equation]{Lemma}
\newtheorem{prop}[equation]{Proposition}
\theoremstyle{definition}
\newtheorem{defn}[equation]{Definition}
\theoremstyle{remark}
\newtheorem{rem}[equation]{Remark}
\theoremstyle{conjecture}

\newcommand{\norm}[1]{\left\Vert#1\right\Vert}
\newcommand{\abs}[1]{\left\vert#1\right\vert}
\newcommand{\set}[1]{\left\{#1\right\}}
\newcommand{\Real}{\mathbb R}
\newcommand{\eps}{\varepsilon}
\newcommand{\To}{\longrightarrow}

\newcommand{\half}{\frac{1}{2}}
\newcommand{\Half}{\dfrac{1}{2}}
\newcommand{\quater}{\frac{1}{4}}
\newcommand{\qedbox}{\hfill \hbox{$\rlap{$\sqcap$}\sqcup$}}
\newcommand{\pf}{\noindent \emph{Proof : }}

\newcommand{\Lie}{\mathrm{Lie}}
\newcommand{\lieso}{\mathfrak{so}}
\newcommand{\lieu}{\mathfrak{u}}
\newcommand{\liesu}{\mathfrak{su}}
\newcommand{\liesp}{\mathfrak{sp}}
\newcommand{\frakg}{\mathfrak{g}}
\newcommand{\frakh}{\mathfrak{h}}
\newcommand{\frakp}{\mathfrak{p}}
\newcommand{\frakk}{\mathfrak{k}}
\newcommand{\frakkp}{\mathfrak{k}^+}
\newcommand{\frakkm}{\mathfrak{k}^-}

\newcommand{\frakt}{\mathfrak{t}}
\newcommand{\fraks}{\mathfrak{s}}
\newcommand{\frakl}{\mathfrak{l}}

\newcommand{\Km}{K^-}
\newcommand{\Kp}{K^+}
\newcommand{\Kpm}{K^\pm}
\newcommand{\wem}{w_-}
\newcommand{\wep}{w_+}
\newcommand{\wepm}{w_\pm}
\newcommand{\lm}{l_-}
\newcommand{\lp}{l_+}
\newcommand{\lpm}{l_\pm}
\newcommand{\Bm}{B_-}
\newcommand{\Bp}{B_+}
\newcommand{\Bpm}{B_\pm}
\newcommand{\ptm}{p_-}
\newcommand{\ptp}{p_+}
\newcommand{\ptpm}{p_\pm}
\newcommand{\qtm}{q_-}
\newcommand{\qtp}{q_+}

\newcommand{\sph}{\mathbb{S}}
\newcommand{\Disk}{\mathbb{D}}

\newcommand{\mul}{\mathrm{mul}}

\newcommand{\Hil}{\mathcal{H}}
\newcommand{\al}{\alpha}

\newcommand{\Tr}{\mathrm{Tr}}
\newcommand{\Ad}{\mathrm{Ad}}
\newcommand{\diag}{\mathrm{diag}}
\newcommand{\id}{\mathrm{id}}
\newcommand{\Id}{\mathrm{Id}}
\newcommand{\diff}{\mathrm{d}}
\newcommand{\Sym}{\mathrm{Sym}}
\newcommand{\Res}{\mathrm{Res}}
\newcommand{\Ind}{\mathrm{Ind}}
\newcommand{\Walp}{W_{\alpha}}
\newcommand{\frakq}{\mathfrak{q}}
\newcommand{\frakqsum}{\mathfrak{q}_1 + \cdots + \mathfrak{q}_\alpha}
\newcommand{\frakno}{\mathfrak{n}_1}
\newcommand{\fraknt}{\mathfrak{n}_2}
\newcommand{\Am}{A_{\_}}
\newcommand{\Ap}{A_{+}}
\newcommand{\Au}{A_u}
\newcommand{\Av}{A_v}

\newcommand{\alfo}{\alpha l+1}
\newcommand{\alft}{\alpha l+2}
\newcommand{\Weyl}{\mathcal{W}}
\newcommand{\Int}{\mathbb{Z}}
\newcommand{\Cpx}{\mathbb{C}}
\newcommand{\Qua}{\mathbb{H}}
\newcommand{\Cay}{\mathbb{O}}
\newcommand{\Cp}{\mathbb{C}\mathrm{P}}
\newcommand{\Hp}{\mathbb{H}\mathrm{P}}

\newcommand{\qi}{\imath}
\newcommand{\qj}{\jmath}
\newcommand{\qk}{\kappa}
\newcommand{\cohom}{\mathrm{cohom}}
\newcommand{\sumi}{\sum_{i=1}^{l-1}}

\newcommand{\fw}{\varpi}

\begin{document}

\title[New Examples of Obstructions]{New Examples of Obstructions to Non-negative Sectional Curvatures
in Cohomogeneity One Manifolds}%
\author{Chenxu He}%
\address{Department of Mathematics, University of Pennsylvania}%
\email{hech@math.upenn.edu}%


\maketitle

\begin{abstract}
K. Grove, L. Verdiani, B. Wilking and W. Ziller gave the first examples of cohomogeneity one
manifolds which do not carry invariant metrics with non-negative sectional curvatures. In this
paper we generalize their results to a larger family. We also classified all class one
representations for a pair $(G,H)$ with $G/H$ some sphere, which are used to construct the
examples.
\end{abstract}

\section{Introduction}
Non-negatively curved Riemannian manifolds have been of interest from the beginning of global
Riemannian geometry. Most examples are obtained via product and quotient constructions, starting
from compact Lie groups. These include all homogeneous spaces and biquotients. Using a gluing
method, J.Cheeger constructed non-negatively curved metrics on the connect sum of any two rank one
symmetric spaces, see \cite{Cheeger}. A breakthrough came with K.Grove and W.Ziller's
generalization of this gluing method to the class of cohomogeneity one manifolds.

A manifold $M$ is called a \emph{cohomogeneity one manifold} if there exists a compact Lie group
$G$ acting on $M$ by isometries and the cohomogeneity of the action, defined as $\cohom(M,G) =
\dim(M/G)$, is equal to $1$. Since the orbit space is one dimensional, it is either a circle or a
closed interval $I$. In the former case, $M$ always carries a $G$ invariant metric with
non-negative sectional curvature. In the latter one, there are precisely two singular orbits $\Bpm$
with isotropy subgroups $\Kpm$ corresponding to the endpoints of $I$, a minimal geodesic between
the singular orbits, and principal orbits corresponding to the interior points with isotropy
subgroup $H$. By the slice theorem, $\Kpm/H$ are spheres $\sph^{\lpm-1}$ ($\lpm \geq 1$), and $M$
can be reconstructed by gluing two disk bundles along a principal orbit as follows
\begin{equation}\label{gluediskbundles}
M = G \times_{\Km}\Disk^{\lm} \cup_{G/H}G \times_{\Kp}\Disk^{\lp},
\end{equation}
where $\Disk^{\lpm}$ is the normal disk to $\Bpm$. Therefore we can identify $M$ with the groups
$H\subset \set{\Km,\Kp} \subset G$ by the gluing construction (\ref{gluediskbundles}).

K.Grove and W.Ziller showed that any cohomogeneity one manifold with codimension two singular
orbits admits a non-negatively curved metric \cite{GZMilnor}. In the same paper, it was conjectured
that any cohomogeneity one manifold admits a non-negatively curved metric. This turns out to be
false. The first examples of an obstruction were discovered by K.Grove, L.Verdiani, B.Wilking and
W.Ziller in \cite{GVWZ}. The most interesting examples are the higher dimensional Kervaire spheres
(of dimension $9$ and up) which are the only exotic spheres that can carry a cohomogeneity one
action, see \cite{St}.

In light of the construction of examples with non-negative sectional curvature and the examples of
obstructions to non-negatively curved metrics, it is important to answer the question \emph{"How
large is the class of cohomogeneity one manifolds that admit a non-negatively curved metric?"}
which was raised by W.Ziller in \cite{Zsurvey1}.

In this paper, we generalize the examples in \cite{GVWZ} to a larger family:
\begin{thmmain} \label{thmmain}
Let $K'/H' = \sph^k$ with $k\geq 2$ and $\rho : K' \To SO(m)$ be a faithful irreducible
representation, which is not the one of $K'$ on $\Real^{k+1}$, and such that $\rho(H')\subset
SO(m-1)$. For any integer $n\geq m+2$, set $G = SO(n)$ and
\begin{eqnarray}
\Km & = & \rho(K')\times SO(n-m) \subset SO(m)\times SO(n-m) \subset SO(n) \nonumber\\
\Kp & = & \rho(H')\times SO(n-m+1)\subset SO(m-1)\times SO(n-m+1) \subset SO(n) \nonumber \\
H & = & \rho(H')\times SO(n-m) \subset SO(n),\nonumber
\end{eqnarray}
then the cohomogeneity one manifold $M$ defined by the groups $H\subset \set{\Km,\Kp}\subset G$
does not admit a $G$ invariant metric with non-negative sectional curvature.
\end{thmmain}

In \cite{GVWZ}, the theorem was proved under the additional assumptions that the slice
representation of $K'$ is not contained in the symmetric square $\Sym^2\rho$ and $\rho(K')$ does
not act transitively on the sphere $\sph^{m-1} = SO(m)/SO(m-1)$.

This theorem is optimal in the sense that if $\rho$ is the representation of $K'$ on $\Real^{k+1}$,
then the manifold $M$ does admit an invariant non-negatively curved metric, since it is
diffeomorphic to the homogeneous space $SO(n+1)/(\rho(K')\times SO(n-m+1))$ endowed with the
cohomogeneity one action of $SO(n) \subset SO(n+1)$. On the other hand, if $n=m+1$ one obtains an
interesting cohomogeneity one manifold where it is not yet known wether that manifold carries a
non-negatively curved invariant metric or not.

The theorem can be extended to the case where $\rho$ is not necessarily irreducible, see Theorem
\ref{thmorth}. A representation $\rho : K'\To SO(m)$ with $\rho(H') \subset SO(m-1)$ is called a
\emph{class one} representation, see Definition \ref{defclassone}. Many class one representations
are not faithful, see Proposition \ref{propclassonetypekernel} and Table \ref{tabletypeclassone},
so they are excluded by the faithfulness requirement in \cite{GVWZ}. However Theorem \ref{thmorth}
allows non-faithful class one representation as a subrepresentation of $\rho$ and it gives us many
more examples, see Section \ref{secgroupcompweylgroup}.

Similar results also hold if $\rho$ is a complex or quaternionic representation, i.e., $G=U(n)$ or
$Sp(n)$ which are stated in the Theorem \ref{thmexamplescomplex} and Theorem
\ref{thmexamplesquater}.

In the above theorem, Theorem \ref{thmorth}, Theorem \ref{thmexamplescomplex} and Theorem
\ref{thmexamplesquater}, the group $K'$ does not need to be connected. But the corresponding groups
$\Kpm$ on the universal cover $\widetilde{M}$ will be connected. Therefore in the rest of the
paper, we can assume that all groups are connected.

\smallskip

The paper is organized as follows. In Section \ref{secprelim}, we briefly recall basic properties
of cohomogeneity one manifolds, the generalized Weyl group and metric properties. Also this section
includes an introduction to class one representations. A more detailed discussion of such
representations is given in Appendix \ref{appclassone}. In Section \ref{secgroupcompweylgroup} we
prove some basic properties of the new examples which will be used in the proof of the orthogonal
case. In Section \ref{secKillingvf}, we study the consequences of the non-negativity assumption on
the metrics. Wilking's rigidity theorem for non-negatively curved manifolds\cite{Wilkingduality}
plays an important role. Section \ref{secprooforth} is devoted to the proof in the orthogonal case.
Using properties of the metrics developed in the previous sections, we draw a contradiction by
looking at sectional curvatures of 2 distinct classes of $2$-planes. In the last section, we sketch
an outline of the proofs in the case where $G=U(n)$ or $Sp(n)$.

\medskip

\emph{Acknowledgements:} The paper is part of the author's Ph.D. thesis at University of
Pennsylvania. The author wants to thank his advisor, Prof. Wolfgang Ziller, for his generous
supports and great patience, and to Prof. Kristopher Tapp for valuable discussions on Wilking's
rigidity results.

\section{Preliminaries}\label{secprelim}

In this section, we recall some basic and well-known facts about cohomogeneity one manifolds. For
more detail, we refer to, for example, \cite{Alekseevsy} and \cite{GWZ}.

As mentioned already, there are precisely two non-principal orbits $\Bpm$ in a simply connected
cohomogeneity one manifold. Suppose $M$ is endowed with an invariant metric $g$ and the distance
between the two non-principal orbits is $L$. Let $c(t)$, $t \in \Real$ be a geodesic minimizing the
distance with $c(0) = \ptm \in \Bm$ and $c(L) = \ptp \in \Bp$. The isotropy subgroups at $\ptpm$
are denoted by $\Kpm$ and the principal isotropy subgroup at any point $c(t)$, $t\in (0,L)$, is
denoted by $H$. We can draw the following group diagram for the manifold $M$:
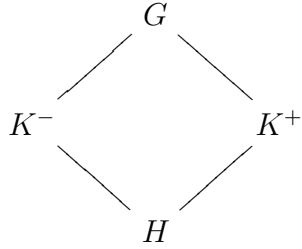
\[
\xymatrix{
                & G             \\
 \Km \ar@{-}[ur]^{ }& &  \Kp \ar@{-}[ul]^{}                 \\
                & H \ar@{-}[ul] \ar@{-}[ur] }
\]
The group diagram $H \subset \set{\Km, \Kp}\subset G$ is not uniquely determined by the manifold
since one can switch $\Km$ with $\Kp$, change $g$ to another invariant metric and choose another
minimal geodesic $c(t)$.

\begin{defn} Two group diagrams are called \emph{equivalent} if they
determine the same cohomogeneity one manifold up to equivariant diffeomorphism.
\end{defn}

The following lemma characterizes which two group diagrams are equivalent, see in \cite{GWZ}.
\begin{lem}\label{LemGroupEquivalent}
Two group diagrams $H \subset \set{\Km, \Kp} \subset G $ and $\tilde{H} \subset \set{\tilde{K}^-,
\tilde{K}^+} \subset G$ are equivalent if and only if after possibly switching the roles of $\Km$
and $\Kp$, the following holds: There exist elements $b \in G$ and $a \in N(H)_0$, where $N(H)_0$
is the identity component of the normalizer of $H$, with $\tilde{K}^- =  b \Km b^{-1}$, $\tilde{H}
= bHb^{-1}$, and $\tilde{K}^+ = ab \Kp b^{-1} a^{-1}$.
\end{lem}

\begin{rem}\label{remequivalentdiag}
If $c(t)$ is the minimal normal geodesic between the two singular orbits, then $b\star c(t)$ is
another minimal geodesic between $\Bpm$ and the associated group diagram is obtained by conjugating
all isotropy groups by the element $b$. We can assume that $b \in N(H) \cap N(\Km)$ in order to fix
$H$ and $\Km$. Conjugation by an element $a$ as in the above lemma usually corresponds to changing
the invariant metric on the manifold.
\end{rem}

Let $C\subset M$ be the image of the minimal geodesic $c(t)$. Then the Weyl group $\Weyl$ of the
$G-$action on $M$ is by definition the stablizer of $C$ modulo its kernel $H$. $\Weyl$ is
characterized in the following proposition.
\begin{prop}\label{propWeyl}
The Weyl group $\Weyl$ of a cohomogeneity one manifold is a dihedral subgroup of $N(H)/H$. It is
generated by involutions $\wepm \in (N(H)\cap \Kpm)/H$ and $C/\Weyl = M/G =[0,L]$. Each of these
involutions can be represented as an element $a\in \Kpm\cap H $ such that $a^2$ but not $a$ lies in
$H$.
\end{prop}

Using the group action, the invariant metric is determined by its restriction to the minimal
geodesic $c(t)$. Suppose $c(t)$ is parameterized by arc length, i.e., $T=\frac{\diff}{\diff t}$ has
length $1$, then we can write $g$ as
\begin{equation*}
g = \diff t^2 + g_t,
\end{equation*}
and $\set{g_t}_{t\in [0,L]}$ is a one-parameter family of homogeneous metrics on the orbits $M_t =
G.c(t)$.

Fix a bi-invariant inner product $Q$ on the Lie algebra $\frakg$ of $G$ and let $\frakp =
\frakh^\perp$ be the orthogonal complement of the Lie algebra $\frakh$ of $H$. For each $X \in
\frakp$, let $X^*$ be the Killing vector field generated by $X$ along $c(t)$, i.e.,$X^*(t) =
\frac{\diff}{\diff s}|_{s=0} \exp(sX).c(t)$. For each $t\in (0,L)$, $M_t$ is diffeomorphic to the
homogeneous space $G/H$, and hence $T_{c(t)}M_t$ can be identified with $\frakp$ by means of
Killing vector fields as $X \mapsto X^*(t)$. Then $g_t$ defines an inner product on $\frakp$ which
is invariant under the isotropy action of $\Ad_H$.  We set
\begin{equation}\label{eqnmetricgt}
g_t(X,Y) = g_t(X^*(t), Y^*(t)) = Q(P_t X, Y) \mbox { for } X, Y \in \frakp,
\end{equation}
where $P_t : \frakp \To \frakp$ is a $Q$-symmetric $\Ad_H$-equivariant endomorphism. The metric $g$
is completely determined by the one parameter family $\set{P_t}$, $t\in [0,L]$, and at $t =0$ and
$L$, $P_t$ should satisfy further conditions to guarantee smoothness of $g$.

On the other hand, each principal orbit $M_t$ is a hypersurface in $M$ with normal vector $T$. If
$S_t X = S_tX^*(t) = -\nabla_{X^*} T $ is the shape operator of $M_t$ at $c(t)$, we have
\begin{equation}\label{eqnshapeop}
S_t = -\Half P_t^{-1}P_t'
\end{equation}
in terms of $P$.

\smallskip

In the rest of this section, we give a short introduction to class one representations with more
details in Appendix \ref{appclassone}. This particular class of representations is well studied,
see, for examples, \cite{VKrepresentation} and \cite{WallachMin}.

First we recall:
\begin{defn}\label{defclassone}
A representation $(\mu, W)$ of a compact Lie group $K$ is called a \emph{real (complex or
quaternionic) representation} if $W$ is a vector space over $\Real$ ($\Cpx$ or $\Qua$).
\end{defn}

\begin{defn}\label{defreptype}
Suppose $W$ is an irreducible representation over the complex numbers, then $W$ is called a
\emph{real} representation or of \emph{real type} if it comes from a representation over reals by
extension of scalars. It is of \emph{quaternionic type} if it comes from a representation over
quaternions by restriction of scalars. It is of \emph{complex type} if it is neither real or
quaternionic.
\end{defn}

For any complex representation $\mu$, let $\mu^*$ denote its complex conjugate. $\mu^*$ is
equivalent to $\mu$ if $\mu$ is of type real or quaternionic and they are non-equivalent if $\mu$
is of complex type.

Suppose $K$ is a compact connected Lie group, then the complexification of a real irreducible
representation $\sigma$ is one of the following
\begin{enumerate}
\item $\sigma_{\Cpx} = \mu$, where $\mu$ is an irreducible representation of real type,
\item $\sigma_{\Cpx} = \mu \oplus \mu^*$, where $\mu$ is irreducible and of complex type,
\label{itemrealcomplex}
\item $\sigma_{\Cpx} = \mu \oplus \mu$, where $\mu$ is irreducible and of quaternionic type.
\label{itemrealquater}
\end{enumerate}
In class (\ref{itemrealcomplex}) and (\ref{itemrealquater}), we often write $\sigma = [\mu]_\Real$.

Conversely, suppose $\mu$ is a complex irreducible representation of $K$ with degree $n$. If $\mu$
is of real type, then there exists a real vector space $\Real^n \subset \Cpx^n$ which is $\mu(K)$
invariant. Let $\sigma$ be the restriction of $\mu$ on $\Real^n$, then $\sigma$ is a real
irreducible representation with $\sigma_{\Cpx} = \mu$. If $\mu$ is not of real type, then we
identify $\Cpx^{n}$ with $\Real^{2n}$ and forget the complex structure on it. Thus $\Real^{2n}$ is
invariant under the $\mu(K)$ action and it gives us a real irreducible representation. We denote
this representation by $\sigma$ and then $\sigma_{\Cpx}$ is either $\mu \oplus \mu^*$ or $\mu
\oplus \mu$ depending on the type of $\mu$.

A quaternionic irreducible representation is obtained by extending the scalar field to the
quaternions for a complex irreducible representation and the converse also holds.

\smallskip

Given an $n$ dimensional real irreducible representation $(\sigma, W)$ of a compact Lie group $K$,
we briefly discuss the classification of the equivariant endomorphisms, i.e., the endomorphism $f:
W \To W$ with $f(\sigma(g).v) = \sigma(g).f(v)$ for any $g\in K$ and $v\in W$.

If $\sigma$ is in class $(1)$, then from Schur's lemma, $f = a \Id$ for some constant $a \in \Real$
and $\Id$ is the identity map of $W$.

If $\sigma$ is in class $(2)$, then $W$ has an orthonormal basis such that $f$ is in the following
form under that basis
\begin{equation}\label{eqnequivarmapcpx}
f=
\begin{pmatrix}
a_0 I_{m} & -a_1 I_{m} \\
a_1 I_{m} & a_0 I_{m}
\end{pmatrix},
\end{equation}
where $a_0$, $a_1$ are constants in $\Real$, $n=2m$ and $I_m$ is the $m\times m$ identity matrix.

If $\sigma$ is in class $(3)$, then there exists an orthonormal basis of $W$ such that $f$ has the
following form under that basis
\begin{equation}\label{eqnequivarmapquater}
f =
\begin{pmatrix}
a_0 I_m & -a_1 I_m & -a_2 I_m & -a_3 I_m \\
a_1 I_m & a_0 I_m & a_3 I_m & -a_2 I_m \\
a_2 I_m & -a_3 I_m & a_0 I_m & a_1 I_m \\
a_3 I_m & a_2 I_m & -a_1 I_m & a_0 I_m
\end{pmatrix},
\end{equation}
where $n=4m$, $a_0$, $a_1$, $a_2$ and $a_3$ are constants in $\Real$.

\smallskip

\begin{defn}
A pair $(K,H)$ of compact Lie groups with $H\subset K$ and $K/H= \sph^k$ ($k\geq 2$) is called a
\emph{spherical pair}.
\end{defn}

If we assume that $K$ is connected, the image of $\mu$ will be a closed subgroup of $SO(l)$, $U(l)$
or $Sp(l)$ if $\mu$ is over $\Real$, $\Cpx$ or $\Qua$.

For each group pair $(K,H)$ with $H \subset K$ a closed subgroup, we have
\begin{defn}
A non-trivial irreducible representation $(\mu, W)$ of $K$ is called \emph{a class one
representation of the pair ($K,H$)} if $\mu(H)$ fixes a nonzero vector $w_0 \in W$.
\end{defn}

\begin{rem}
From Proposition \ref{propclassonetypekernel}, the class one representations of spherical pairs are
almost effective, i.e., the kernel of $\mu$ is discrete, except for the following cases:
\begin{enumerate}
\item $(K,H)=(U(n),U(n-1)_m)$ and $\mu = ae_1 - ae_n$ where $a$ is a positive integer. The kernel is the diagonal
embedded $U(1) \subset U(n)$.
\item $(Sp(n)\times Sp(1), Sp(n-1)\times Sp(1))$ and $\mu = a \fw_2$ where $a$ is a positive
integer. The kernel contains the $Sp(1)$ factor.
\item $(Sp(n)\times U(1), Sp(n-1)\times U(1)_m)$ and $\mu = (a\fw_1 + b\fw_2) \otimes \id$ where
$a$, $b$ are non-negative integers with $a+b\geq 1$ and $\id$ is the trivial representation of
$U(1)$. The kernel contains the $U(1)$ factor.
\end{enumerate}

More information on the non-trivial kernels for class one representations is given in Proposition
\ref{propclassonetypekernel} and listed in Table \ref{tabletypeclassone}.
\end{rem}

In the case where $(K,H)= (SO(k+1),SO(k))$, let $\fw_1$ be the highest weight of the standard
representation $\varrho_{k+1}$ on $\Real^{k+1}$, then the class one representations over $\Real$
are precisely those with the highest weights as $m\fw_1$, $m=1, 2, \ldots$. These representation
spaces can be realized as the space of homogeneous harmonic polynomials. Let $\set{x_1, \ldots,
x_{k+1}}$ be the basis of $\Real^{k+1}$ and $SO(k+1)$ act by the matrix multiplication. Then an
element $A \in SO(k+1)$ acts on a polynomial $f(x_1, \ldots, x_{k+1})$ through the action on the
variables, i.e., $(A.f)(x_1, \ldots, x_{k+1}) = f(A^{-1}.x_1, \ldots, A^{-1}.x_{k+1})$.

A polynomial $f$ is called harmonic if $\Delta f = 0$ where $\Delta = \sum_{i=1}^{k+1}
\frac{\diff^2}{\diff x_i^2}$. Let $H_m$ be the space of homogeneous harmonic polynomials in $x_1,
\ldots, x_{k+1}$ of degree $m$, then the representation of $SO(k+1)$ on $H_m$ has the highest
weight $m\fw_1$. All of them are of real type. If $k+1$ is odd, then the representation for any
positive $m$ is faithful. If $k+1$ is even, then the class one representation is faithful if and
only if $m$ is odd. As we will see that the class one representations of other spherical pairs
$(K,H)$ are the irreducible components when $H_m$ is restricted to the subgroup $K$ of $SO(k+1)$.

The real, complex and quaternionic class one representations of the spherical pairs are classified
in Appendix \ref{appclassone}. Further properties are discussed there as well.

\section{Weyl group and smoothness}\label{secgroupcompweylgroup}

First let us state the most general result in the orthogonal case as a generalization of the
theorem in the introduction using the concept of class one representation:
\begin{thm}\label{thmorth}
Let $K'/H' = \sph^k$ with $k\geq 2$ and $\rho : K' \To SO(m)$ be a faithful representation that
contains a (not necessarily faithful) irreducible class one representation $\mu$ of the pair
$(K',H')$ such that one of the followings holds:
\begin{enumerate}
\item $\deg \mu \geq k+2$ if $\mu$ is of real type or the multiplicity of $\mu$ in $\rho$,
denoted by $\mul(\mu,\rho)$, is equal to $1$;
\item $\deg \mu \geq 2(k+2)$ if $\mu$ is of complex type and $\mul(\mu,\rho)\geq 2$;
\item $\deg \mu \geq 4(k+2)$ if $\mu$ is of quaternionic type and $\mul(\mu,\rho) \geq 2$.
\end{enumerate}

We assume that $n \geq m+2$ if $\mul(\mu,\rho) = 1$, and that $n \geq m+3$ if $\mul(\mu,\rho) \geq
2$. If we set $G = SO(n)$ and
\begin{eqnarray}
\Km & = & \rho(K')\times SO(n-m) \subset SO(m)\times SO(n-m) \subset SO(n) \nonumber\\
\Kp & = & \rho(H')\times SO(n-m+1)\subset SO(m-1)\times SO(n-m+1) \subset SO(n)
\label{groupdiagram}\\
H & = & \rho(H')\times SO(n-m) \subset SO(n),\nonumber
\end{eqnarray}
then the cohomogeneity one manifold $M$ defined by the groups $H\subset \set{\Km,\Kp}\subset G$
does not admit a $G$ invariant metric with non-negative sectional curvature.
\end{thm}

\begin{rem}
Since we do not assume that $\mu$ is faithful, a non-faithful class one representation is allowed
in this construction. In other words, if $\mu$ is a class one representation with $\deg \mu \geq
k+2$, we choose a representation $\tau$ with $\ker \tau \cap \ker \mu =\set{1}$. Then $\rho = \tau
\oplus \mu$ satisfies the conditions in the theorem. For example, take $(K',H') = (SO(6),SO(5))$
with $k = 5$ and let $\mu$ be the $20$ dimensional representation of $SO(6)$ with the highest
weight $2\fw_1$. To construct a cohomogeneity one manifold, we can choose $\tau = \varrho_6$ as the
standard representation of $SO(6)$ which is faithful and let $\rho = \tau \oplus \mu$.
\end{rem}

\begin{rem}
If $\rho$ contains only one copy of $\mu$ or $\mu$ is of real type, then Table
\ref{tableclassonedimcomparereal} in Proposition \ref{propdimcomparereal} lists all real class one
representations which have dimensions smaller than $k+2$. We see that only the following
representations are excluded by the assumption, $\deg \mu \geq k+2$: the defining representation of
$K'$ on $\Real^{k+1}$, the $9$ dimensional representation $\varrho_9$ of the pair ($Spin(9)$,
$Spin(7)$), the $5$ dimensional representation of $(Sp(2),Sp(1))$ and the $3$ dimensional
representations of $(SU(2), \set{\Id})$ and $(U(2),U(1))$. All of these representations are not
faithful, so one needs to add another representation $\tau$ with $\ker \tau \cap \ker \mu
=\set{\Id}$ to define a cohomogeneity one manifold. For such manifolds, we do not know if they
admit an invariant metric with non-negative curvature.
\end{rem}

\begin{rem}
If $\rho$ contains more than one copy of $\mu$ and $\mu$ is not of real type, then the further
restriction $\deg \mu \geq 2 (k+2)$ or $4 (k+2)$ excludes $8$ more representations as listed in the
last part of Table \ref{tableclassonedimcomparereal}. Among them, the following representations are
faithful: $\mu = [2\fw_1]_\Real$ for the pair $(SU(3),SU(2))$, $\mu=[\fw_1 +\fw_2]_\Real$ for the
pair $(Sp(2),Sp(1))$ and $\mu = [3\fw_1]_\Real$, $[5\fw_1]_\Real$ and $[7\fw_1]_\Real$ for the pair
$(Sp(1),\set{1})$. The first one is of complex type and the other four are of quaternionic type.
They can be used to construct cohomogeneity one manifolds without adding other representations, for
example, take $(K',H') = (Sp(1),\set{1})$ and $\rho = [3\fw_1]_\Real \oplus [3\fw_1]_\Real$. Thus
Theorem \ref{thmorth} does not give obstruction for such manifolds.
\end{rem}

\begin{rem}
The lowest dimensional example of Theorem \ref{thmorth} is obtained as follows. Take $(K',H')=
(SO(3),SO(2))$ with $K'/H' = \sph^2$, then the lowest dimensional class one representation $\mu$
with $\deg > 3$ is the unique $5-$dimensional representation of $SO(3)$ which is also faithful. If
$m=5$ and $n=7$, the manifold $M$ has dimension $20$ and isotropy groups
\begin{equation*}
\mu(SO(2))\times SO(2) \subset \set{\mu(SO(3))\times SO(2), \mu(SO(2)) \times SO(3)} \subset SO(7).
\end{equation*}
Notice that this example is already covered by Theorem 3.2 in \cite{GVWZ}.
\end{rem}

\smallskip

We now describe the explicit embedding of the groups. $SO(m)\times SO(n-m)$ is embedded in
$G=SO(n)$ block-wise, i.e., $SO(m)$ sits in the upper-left $m\times m$-block and $SO(n-m)$ is in
the lower-right block. By assumption, $\rho$ is a faithful orthogonal representation of $K'$ with
representation space $V=\Real^m$. Let $W_1, \ldots, \Walp$ be invariant subspaces of $V$ such that
they are pairwisely orthogonal, and the restrictions of $\rho(K')|_{W_i} = \mu$ are equivalent and
irreducible. Let $U$ be the orthogonal complement of $W = W_1 \oplus \cdots \oplus \Walp$.
According to the decomposition of $V$ into invariant spaces, we can write $\rho$ as
\begin{equation}\label{eqndecomprep}
\rho = \tau \oplus \mu \oplus \cdots \oplus \mu,
\end{equation}
where $\tau$ is restriction of $\rho$ to $U$ and $\mu$ is the class one representation of
$(K',H')$. Furthermore, $\mu$ is not a subrepresentation of $\tau$. Suppose $r=\dim U$, $l=\dim W_i
= \deg \mu$, then $m=r+ \alpha l$. By choosing a suitable basis of $V$, $\rho(x)$ for $x\in K'$, is
a block diagonal matrix in $SO(m)$:
\begin{equation}\label{eqnrhox}
\rho(x) =
\begin{pmatrix}
\tau(x) &         &        &        \\
        & \mu(x)  &        &        \\
        &         & \ddots &        \\
        &         &        & \mu(x) \\
\end{pmatrix}
\in SO(r) \times SO(l) \times \cdots \times SO(l).
\end{equation}

When $\mu$ is restricted to the subgroup $H'\subset K'$, it is not irreducible any more. Let $m_0$
be the multiplicity of the trivial representation in the restriction $\mu |_{H'}$. Hence for any
element $y \in H'$, $\mu(y) \in SO(l-m_0) \subset SO(l)$ and $SO(l-m_0)$ is embedded as the upper
left $(l-m_0)\times (l-m_0)$ block in $SO(l)$.

\smallskip

With the explicit description of the embeddings, we can prove the following proposition on the Weyl
group of the new examples.
\begin{prop}\label{propweyl}
The Weyl group $\Weyl$ is isomorphic to $\Int_2 \times \Int_2 $.
\end{prop}

\pf For any element $x\in K'$ and $A\in SO(n-m)$, let $M(x,A)$ denote the block diagonal matrix
$\diag(\rho(x), A)$ with $\rho(x) \in SO(m)$. If $x \in H'$, $A$ can be considered as a matrix in
$SO(n-m+1)$ since $\rho(x) \in SO(m-1)$.

First notice that $\wep$ can be represented by the an element $a \in \Kp = \rho(H')\times
SO(n-m+1)$ which is not in $H$, but $a^2 \in H$. Let
\begin{equation}\label{eqnwep}
\wep = M(\id, \left(\begin{smallmatrix}
-1 &    & \\
   & -1 & \\
   &    & I_{n-m-1}
\end{smallmatrix}\right)) = a=
\begin{pmatrix}
 I_{m-1} &     &        &        &     \\
         & -1  &        &        &     \\
         &     & -1 &        &     \\
         &     &        & I_{n-m-1} &     \\
\end{pmatrix},
\end{equation}
where $I_{k}$ is the $k\times k$ identity matrix.

Suppose $b = M(x,A)$ is a representative of $\wem$, i.e., $b \in N_{\Km}(H)$ and $b^2$ but not $b
\in H$. In the following we will determine the element $b$ in three different cases depending on
the class one representation $\mu$.

\textsc{Case 1:} $\mu$ is of real type and $m_0 = 1$. In each $W=W_i$($i =1, \ldots, \alpha$), let
$v$ be a unit vector fixed by $\mu(H')$(or equivalently by $H$) and $X$ be its orthogonal
complement. Then $\dim X = l-1$. Since $b \in N_{\Km}(H)$, $b.v$ is also fixed by $H$ and has the
same length as $v$, i.e., $b.v= \pm v$. Since $b$ is an orthogonal transformation, $b$ maps $X$ to
itself, i.e., when restricted on $W$, $b=\left(\begin{smallmatrix} b_1 & 0
 \\ 0 & \det b_1  \end{smallmatrix}\right)$ where $b_1 \in O(l-1)$. Therefore the representative
$b$ has the following matrix form:
\begin{equation}\label{eqnwem}
\wem = b =
\begin{pmatrix}
A_1 &     & & & & &\\
    & A_2 & & & & &\\
    &     & \det A_2 & & & &\\
    &     &          & \ddots &  &  &\\
    &     &          &        & A_2 & & \\
    &     &          &        &     & \det A_2 & \\
    &     &          &        &     &          & I_{n-m}
\end{pmatrix},
\end{equation}
where $A_1 = \tau(x) \in SO(r)$, $A_2 = b_1 \in O(l-1)$ and there are $\alpha$ copies of
$\left(\begin{smallmatrix}
A_2 & \\
& \det A_2
\end{smallmatrix}\right)$.

\textsc{Case 2:} $\mu$ is not of real type and $m_0 = 2$. In this case we have $\mu = \nu \oplus
\nu^*$ where $\nu$ is a complex class one representation for the pair $(K',H')$. If $\nu$ is of
quaternionic type, then $\nu^* = \nu$ and $\mu = \nu \oplus \nu$. As in Case 1, we have $W= X
\oplus_\perp Y$ where $Y$ is the $2$ dimensional subspace fixed by $\mu(H')$ and $X$ is its
orthogonal complement. The orthogonal transformation $\mu(x)$ maps $X$ and $Y$ to themselves and
the matrix $\mu(x)$ has the form
$\left(\begin{smallmatrix} b_1 & 0 \\
0 & b_2 \end{smallmatrix}\right) \in S(O(l-2)\times O(2))$. Since $\mu(x)^2 \in H$ implies $b_2^2 =
I_2$, $b_2$ is a symmetric matrix. Since $\mu = \nu \oplus \nu^*$, $b_2$ commutes with the matrix
$\left(\begin{smallmatrix} 0 & -1 \\
1 & 0 \end{smallmatrix}\right)$ which implies $b_2 = \pm I_2$. Therefore $b$ also has the matrix
form as in (\ref{eqnwem}).

\textsc{Case 3: }$\mu$ is not in Case 1 or Case 2. From the classification of class one
representations in Theorem \ref{thmclassonecpx} and their types in Proposition
\ref{propclassonetypekernel}, the only representation of this type is $(K',H') = (Sp(k),Sp(k-1))$
where $\mu$ has the highest weight $p \fw_1 + q \fw_2$ with $p\geq 1$. Here $\fw_1$ and $\fw_2$ are
the 1st and 2nd fundamental weights of $Sp(n)$. Since $\rho$ is a faithful representation, we have
$x\in K'\cap N_{K'}(H') - H'$ and $x^2 \in H'$. Take $x = -I_k \in Sp(k)$, then $x$ satisfies these
restrictions. If we view $Sp(k)$ as a subgroup of $SO(4k)$, then $\mu$ is contained in the
restriction of some class one representation $\nu$ of $SO(4k)$ to $Sp(k)$ by Theorem
\ref{thmharmonic}. The representation space of $\nu$ is consisted of homogeneous harmonic
polynomials, so the image $\nu(x)$ is $\pm \Id$ depending on the parity of the degree of the
polynomials. Therefore $\mu(x)$ is equal to $\pm \Id$ and the element $\wem$ can be represented by
the following matrix:
\begin{equation}\label{eqnwemsp}
\wem = b = \begin{pmatrix}
A_1 &     & &\\
    & \eps I_{\alpha l} & &\\
    &     &          & I_{n-m}
\end{pmatrix},
\end{equation}
where $A_1 = \tau(x) \in SO(r)$ and $\mu(x) = \eps I_l$ with $\eps = \pm 1$.

In each case, from the given representatives of $\wepm$, it is easy to check that $\wep \wem = \wem
\wep$ which is not an element in $H$. Thus $\Weyl =<\wem, \wep> \cong \Int_2 \times \Int_2$.
\qedbox

\begin{rem}
Let $a$ be an element in $N(H)_0$ which does not lie in $N(\Km)$ or $N(\Kp)$ and let $\bar{M}$ be
the cohomogeneity one manifold defined by the group diagram $H \subset \set{\Km, a\Kp a^{-1}}
\subset G$. As pointed out in Remark \ref{remequivalentdiag}, $M$ and $\bar{M}$ usually have
different Weyl groups and different invariant metrics though they are $G$-equivariantly
diffeomorphic. Therefore, if one family of invariant metrics does not admit non-negative curvature,
it does not necessarily follow that the other family is obstructed as well. In our example, if
 the multiplicity of the trivial representation in $\mu|_{H'}$ is equal to one, i.e.,
$m_0 = 1$, then our arguments which show the obstructions work for all equivalent diagrams. This is
the case for most of the class one representations of spherical pairs, see, for example, Theorem
\ref{thmclassonecpx}. On the other hand, if $m_0 \geq 2$, then we have to put the further
restriction on the diagram: the Lie algebra of $\Kp$ contains the subspace
$\mathrm{span}\set{E_{m,m+1}, \ldots, E_{m,n}}$. Here we use $E_{i,j}$, $i\neq j$, to denote the
skew-symmetric matrix having $1$ in the $i,j$-entry, $-1$ in the $j,i$-entry and zero otherwise.

\end{rem}

\smallskip

Using the explicit representatives of the generators of the Weyl group, we have the following
smoothness condition of an invariant metric on $M$.
\begin{lem}\label{lemweylsymmetry}
For any $G-$invariant metric $g$ on $M$, let $h(t)$ be the length of the Killing vector field
generated by $E_{m,m+1}$ along the geodesic $c(t)$. Then $h(t)$ is an even function with $h(0)\neq
0$ and $h(L) = 0$.
\end{lem}
\pf The fact that $E_{m,m+1}$ lies in the Lie algebra of $\Kp$ but not $\Km$ implies that $h(0)\neq
0$ and $h(L) = 0$. The generator $\wem$ is a reflection of $c(t)$ at the point $\ptm$ and maps
$c(t)$ to $c(-t)$. The induced map $\diff \wem$ takes $T_{c(t)}M$ to the tangent space
$T_{c(-t)}M$. From the matrix form (\ref{eqnwem}) and (\ref{eqnwemsp}) of a representative of
$\wem$, we have $\diff \wem (E^*_{m,m+1}(t))$ = $\pm E^*_{m,m+1}(-t)$. Therefore $h(t)=h(-t)$,
i.e., $h(t)$ is an even function. \qedbox

\smallskip

\smallskip

\section{Restrictions on the Metric Along the Normal Geodesic}\label{secKillingvf}

In this section, we begin the study of the invariant metrics in our examples. In general, the
family of $G$-invariant metrics on $M$ is very large. There are many rigidity results for
positively curved metrics. For example, in even dimensions, L. Verdiani classified all positively
curved cohomogeneity one manifolds, see \cite{VerdianiEven1} and \cite{VerdianiEven2}. In odd
dimensions, K.Grove, B.Wilking and W.Ziller obtained a short list of cohomogeneity one manifolds
which possibly have invariant positively curved metrics in \cite{GVWZ}. Recently K.Grove,
L.Verdiani and W.Ziller have succeeded in constructing positively metric on one of them in
\cite{GVZPos}. On the other hand, there are few rigidity results for non-negatively curved metrics,
even in the cohomogeneity one case. Recently, B.Wilking proved some fundamental rigidity theorems
for non-negatively curved manifolds in a general setting which will play an important role in our
proof.

\smallskip

Let $c: \Real \To M$ be a geodesic and let $\Lambda$ be an $(n-1)$-dimensional family of normal
Jacobi fields. The Ricatti operator $L(t)$ is the endmorphism of the normal bundle
$(\dot{c}(t))^\perp$ defined by $L(t)J(t) = J'(t)$ for $J \in \Lambda$. From Theorem 9 and
Corollary 10 in \cite{Wilkingduality}, We have

\begin{thm}[Wilking's Rigidity Theorem]\label{thmwilking}
Suppose the Riccati operator for an $(n-1)$-dimensional family $\Lambda$ of normal Jacobi fields is
self-adjoint, and define the smooth subbundle $\Upsilon$ of $(\dot{c}(t))^\perp$ by:
\[
 \Upsilon =
\mathrm{span}\set{J \in \Lambda | J(t) = 0 \mbox{ for some }t \in \Real}.
\]
Then, if $M$ has non-negative curvature, we have:
\begin{equation}
\Lambda = \Upsilon \oplus \set{J \in \Lambda | J \mbox{ is parallel }}, \label{eqndecompJacobi}
\end{equation}
and
\[
\Upsilon(t) = \set{J(t) |J \in \Upsilon}\oplus\set{J'(t)|J \in \Upsilon, J(t) =0}.
\]
A point $t_0\in \Real$ or $c(t_0)$ is said to be \emph{singular} if $J(t_0) = 0$ for some $J\in
\Upsilon$. Otherwise $t_0$ is said to be \emph{generic}. Thus if $J\in \Lambda$ and $J(t_0) \perp
\Upsilon(t_0)$ at a generic $t_0$, then $J$ is parallel along $c(t)$, $t\in \Real$.
\end{thm}

\begin{rem}
If there exists a subbundle $E \subset (\dot{c}(t))^\perp$ that is invariant under parallel
transport, then Theorem \ref{thmwilking} can also be applied to a $\mathrm{rank}$ $E$ dimensional
family of normal Jacobi fields in $E$ with self-adjoint Ricatti operator.
\end{rem}

In our example, let $c(t)$ be the minimal geodesic between $\Bpm$, and $X^*(c(t))$, $X\in
\frakh^\perp$, an $(n-1)$ dimensional family of Jacobi fields. Its Riccati operator $L(t)$ is self
adjoint, since it is equal to the shape operator $- \half P_t^{-1}P'_t$, see (\ref{eqnshapeop}). We
will apply Theorem \ref{thmwilking} to a subfamily of these Jacobi fields.

Let $\frakg = \lieso(n)$ and $\frakh$ be the Lie algebras of $G=SO(n)$ and $H=\rho(K')\times
SO(n-m)$ respectively. Choose the bi-invariant inner product $Q= -\half \Tr$ on $\frakg$ for which
$\set{E_{i,j}}$ is an orthonormal basis and let $\frakp$ be the orthogonal complement of $\frakh
\subset \frakg$. First we identify some subspaces of $\frakp$.

Let
\begin{eqnarray}
\frakq_0 & = & \mathrm{span} \set{E_{i,j}|1\leq i\leq r, m+1\leq j \leq n }\nonumber \\
\frakq_1 & = & \mathrm{span} \set{E_{i,j}|r+1\leq i\leq r+l, m+1\leq j \leq n }\label{eqnsubspaceq} \\
         & \cdots & \nonumber\\
\frakq_\alpha & = & \mathrm{span} \set{E_{i,j}|m+1-l\leq i\leq m, m+1\leq j \leq n }\nonumber,
\end{eqnarray}
and $\frakq = \frakq_0 + \frakqsum$. We write the last subspace $\frakq_\alpha$ as a sum of two
subspaces as follows:
\begin{eqnarray}
\frakno & = & \mathrm{span}\set{E_{i,j}| m+1-l\leq i \leq m -1, m+1 \leq j \leq n},\label{eqnsubspacen}\\
\fraknt & = & \mathrm{span}\set{E_{m,j}| m +1 \leq j \leq n}\nonumber .
\end{eqnarray}

Let $\frakq^\perp$ be the $Q$-orthogonal complement of $\frakq$ in $\frakp$. Since $\frakq^\perp$
is the fixed point set by the isotropy action of the subgroup $SO(n-m) \subset H = \rho(H')\times
SO(n-m)$, Schur's lemma implies that the Killing vector field $X^*$, $X \in \frakq^\perp$, is
orthogonal along $c(t)$ to $Y^*$ for any $Y \in \frakq$.

\smallskip

\noindent \emph{Terminology. } In the rest of the paper, for any two subspaces $\frakp_1$,
$\frakp_2 \subset \frakp$, the notation $\frakp_1^* \perp \frakp_2^*$ means that any Killing vector
field generated by an element in $\frakp_1$ is orthogonal to any Killing vector field generated by
an element in $\frakp_2$ along $c(t)$.

\smallskip

Since parallel translation commutes with the action of $\Ad_{H}$, and since $\frakq^\perp$ is the
fixed point set of $SO(n-m)$ in $H$, it follows that $(\frakq^\perp)^*$, and hence also $\frakq^*$,
is invariant under parallel translation. By the same reasoning, $P^{-1}P'$ preserves $\frakq^*$ and
thus $\frakq^*$ forms a self adjoint family of Jacobi fields to which we can thus apply Theorem
\ref{thmwilking}.

We determine the component $\Upsilon$ in the splitting (\ref{eqndecompJacobi}) of the Jacobi fields
$\frakq^*$. $\ptp = c(L)$ is a singular point. The element $\wep$ fixes $\ptp$ and reflects $c(t)$
about $\ptp$. Let $\qtm = c(2L) = \wep(\ptm) \in \Bm$, then the isotropy subgroup at $\qtm$ is
$\Km_1 = \Ad_{\wep} \Km$ with Lie algebra $\frakkm_1 = \Ad_{\wep} \frakkm$. Similarly, $\wem$ fixes
$\ptm$ and reflects $c(t)$ about $\ptm$. Let $\qtp = c(-L) = \wem(\ptp) \in \Bp$, then $\qtp$ has
isotropy subgroup $\Kp_1 = \Ad_{\wem} \Kp$ with Lie algebra $\frakkp_1 = \Ad_{\wem} \frakkp$. Since
$\wem.\wep = \wep.\wem$, the image of $\qtm$ under the reflection $\wem$ about $\ptm$ is $\wem
(\qtm) = \wem. \wep(\ptm) = \wep.\wem(\ptm) = \wep(\ptm) = \qtm$, i.e., $c(2L) = c(-2L)$. Therefore
$c(t)$ is a closed geodesic with period $4L$ and the singular points are $\ptp = c(L)$ and $\qtp =
c(3L)$. The vanishing Killing vector fields are those generated by the vectors in the Lie algebras
of the isotropy subgroups at singular points. Notice that if $X \in \frakq$, then $X^*(\ptm) \ne 0$
and $X^*(\qtm) \ne 0$. Theorem \ref{thmwilking} implies

\begin{lem}\label{lemvanshingKillingVF}
If $Y \in \frakq$ such that $Y^*(\ptm) \perp X^*(\ptm)$ for all $X \in \fraknt$, then $Y^*$ is a
parallel Jacobi field along $c(t)$.
\end{lem}

In the following, we prove some properties of the invariant metrics $g$ on $M$ under the
non-negative sectional curvatures assumption.

\smallskip

First we observe

\begin{lem}\label{lemq0orth}
Suppose $(M,g)$ is non-negatively curved, then $\frakq_0^*$ is orthogonal to $(\frakqsum)^*$ along
$c(t)$.
\end{lem}

\pf At the generic point $\ptm=c(0)$, the metric $g$ restricted to the singular orbit $\Bm \cong
G/\Km$ is $\Ad_{\Km}$ invariant. The actions of $\Ad_{\Km}$ on $\frakq_0$ and $\frakq_i$ $(i>0)$,
are $\tau \otimes \rho_{n-m}$ and $\mu \otimes \rho_{n-m}$ respectively, where $\varrho_{n-m}$ is
the standard representation of $SO(n-m)$ on $\Real^{n-m}$. Since $\tau$ does not contain $\mu$ as a
subrepresentation, $\frakq_0^*$ is orthogonal to $\frakq_i^*$ at $\ptm$. In particular $\frakq_0^*$
is orthogonal to $\fraknt^*$, so any Killing vector field generated by a vector in $\frakq_0$ is
parallel along $c(t)$. Hence $\frakq_0^*$ is orthogonal to $(\frakqsum)^*$ along $c(t)$. \qedbox

\smallskip

Next we study the metric $g$ on the space $(\frakqsum)^*$. Recall $h^2(t) =
g(E^*_{m,m+1},E^*_{m,m+1})$ is an even function with $h(0) \ne 0$ and W.L.O.G. we may assume that
$h(0) = 1$.

\subsection{$\mu$ is of real type or the multiplicity of $\mu$ in $\rho$ is one.}

We first consider the case where $\mu$ is of real type. We denote $E_{r+(i-1)l+a,m+\xi}$ by
$E_{a,i,\xi}$ for $a = 1, \ldots, l$, $i =1, \ldots, \alpha$ and $\xi = 1, \ldots, n-m$.

Since $\Ad_{\Km}$ commutes with $P(0)$, the restriction of $P(0)$ to $\frakq_i$, composed with the
projection to $\frakq_j$, is an equivalence between the $\Km$ irreducible representations
$\frakq_i$ and $\frakq_j$. Since they are orthogonal, Schur's Lemma implies that $P(0)(E_{a,i,\xi})
= f_{i,j} E_{a,j,\xi}$ for some constant $f_{i,j}\in \Real$. Furthermore, $f_{i,j}=f_{j,i}$ from
the $Q$-symmetry of $P_t$. In terms of inner product of Killing vector fields, we have
\begin{equation}\label{eqnfijreal}
f_{i,j} = g(E_{1,i,1}^*, E_{1,j,1}^*)_{c(0)}.
\end{equation}
The assumption $h(0) = 1$ implies that $f_{\alpha,\alpha} = 1$.

\begin{lem}\label{lemorth}
Suppose $(M,g)$ is non-negatively curved, then we have
\begin{enumerate}
\item $E_{a,i,\xi}^*$ is a parallel Jacobi field along $c(t)$ if $a\ne l$;
\item $E_{a,i,\xi}^*$ is orthogonal to $E_{b,j,\zeta}^*$ along
$c(t)$ if $a\neq b$ or $\xi \neq \zeta$;
\item $E_{a,i,\xi}^*$ has the same length as $E_{a,i,\zeta}^*$
along $c(t)$;
\item At the point $\ptm = c(0)$, $P_0(E_{a,i,\xi})= \sum_{j=1}^{\alpha} E_{a,j,\xi} f_{i,j}$.
\end{enumerate}
\end{lem}

\pf At the generic point $c(0)$, the $\Ad_{\Km}$ actions on $\frakq_i$ and $\frakq_j$ are
equivalent and given by the irreducible representation $\mu \otimes \varrho_{n-m}$, from Schur's
lemma and the fact that $\mu$ is of real type or the multiplicity of $\mu$ in $\rho$ is one,
$E^*_{a,i,\xi}(0)$ is orthogonal to $\Upsilon(0)= \mbox{span} \set{E^*_{l,\alpha,\varsigma}(0)|
\varsigma = 1,\ldots, n-m}$ for $a\ne l$. Hence it is a parallel vector field from the last part of
Theorem \ref{thmwilking} which proves (1).

On each principal orbit $M_t \cong G/H$, $\Ad_H$ acts on each $\frakq_i$ ($i>0$), by the
representation $Res^{K'}_{H'}(\mu) \otimes \varrho_{n-m}$. By Schur's lemma we have $E^*_{a,i,\xi}$
is orthogonal to $E^*_{b,j,\zeta}$ along $c(t)$ if $\xi \neq \zeta$ and $E^*_{a,i,\xi}$ has the
same length as $E^*_{a,i,\zeta}$. This proves (3) and one case of (2) where $\xi \neq \zeta$.

Suppose $\zeta = \xi$ and $a\neq b$. If none of $a$ or $b$ is equal to $l$, then the two vector
fields $E_{a,i,\xi}^*$ and $E^*_{b,j,\xi}$ are parallel from (1). Using Schur's lemma again and the
fact that $a\ne b$, they are orthogonal to each other at $c(0)$ and then along the normal geodesic
$c(t)$.

If one of $a$ and $b$, say $b$, is equal to $l$, then $E^*_{a,i,\xi}$ is a parallel vector field.
Write $E_{l,j,\xi}^*(0) = (E_{l,j,\xi}-\lambda E_{l,\alpha,\xi})^*(0) + \lambda E_{l,\alpha,
\xi}^*(0)$, where the constant $\lambda$ is determined by the following equation:
\begin{equation*}
g(E_{l,j,\xi}^*, E_{l,\alpha,\xi}^*)_{c(0)} = \lambda g(E_{l,\alpha,\xi}^*,
E_{l,\alpha,\xi}^*)_{c(0)}.
\end{equation*}
Thus $(E_{l,j,\xi}-\lambda E_{l,\alpha,\xi})^*(0) \perp \Upsilon(0)$ and hence
$(E_{l,j,\xi}-\lambda E_{l,\alpha,\xi})^*$ is a parallel vector field. Furthermore, $E_{a,i,\xi}^*$
is orthogonal to $(E_{l,j,\xi}-\lambda E_{l,\alpha,\xi})^*$ at $c(0)$, so they are orthogonal to
each other along $c(t)$. Thus $E_{a,i,\xi}^*$ is orthogonal to $E^*_{l,j,\xi}$ along $c(t)$.

The formula of $P_0(E_{a,i,\xi})$ in $(4)$. follows easily from the defining equation
(\ref{eqnfijreal}) of $f_{i,j}$ and $(2)$. \qedbox

\medskip

From Lemma \ref{lemorth} above, the restriction of the endomorphism $P$ on $(\frakqsum)^*$ at $t=0$
has the following matrix form:
\begin{equation}\label{eqnPinitial}
P_0 =
\begin{pmatrix}
f_{1,1} I_l & \cdots & f_{1,\alpha}I_l\\
\vdots & \ddots & \vdots\\
f_{\alpha,1}I_l & \cdots & f_{\alpha, \alpha}I_l
\end{pmatrix}.
\end{equation}

We have seen that there are plenty of parallel Killing vector fields in $(\frakqsum)^*$. Using
these parallel vector fields, we can determine the restriction of $P_t$ on $(\frakqsum)^*$.

\begin{thm}\label{thmPtrealtype} If the cohomogeneity one manifold $(M,g)$ has non-negative sectional
curvature and the class one representation $\mu$ is of real type, then for any $i = 1, \ldots,
\alpha$ and $\xi = 1, \ldots, n-m$, we have
\begin{enumerate}
\item $P_t(E_{a,i,\xi}) = \sum_{j=1}^{\alpha} f_{i,j} E_{a,j,\xi}$, for $a = 1, \cdots, l-1$;
\item $P_t(E_{l,i,\xi}) = \sum_{j=1}^{\alpha}p_{i,j}(t) E_{l,j,\xi}$ and $p_{i,j}(t)$ is defined as
\begin{equation}\label{eqnpijreal}
p_{i,j}(t) = (h^2(t)-1) a_i a_j  + f_{i,j},
\end{equation}
where $a_i = f_{i,\alpha}$ and $a_\alpha = 1$.
\end{enumerate}
\end{thm}

\pf Part (1) is obvious since every component of both $E^*_{a,i,\xi}$ and $E^*_{a,j,\xi}$ are
parallel vector fields along $c(t)$ if $a\leq l-1$.

For part (2), let
\begin{equation*}
X_i = E_{l,i,\xi} - a_i E_{l,\alpha,\xi}, \quad i = 1, ..., \alpha.
\end{equation*}
Then $X_i \in \frakqsum$ and generates a Killing vector field $X_i^*$ along $c(t)$. By the
definition $a_i = f_{i,\al}$, the defining equation of $f_{i,\alpha}$ in (\ref{eqnfijreal}) and
$(4)$ in Lemma \ref{lemorth}, we have
\begin{equation*} g(X_i^*, E_{l,\alpha, \xi}^*)_{c(0)} = 0,
\end{equation*}
or $X_i^*(0) \perp \Upsilon(0)$. Therefore $X_i^*$ is a parallel vector field. By the formula
(\ref{eqnshapeop}) of the shape operator, we have
\begin{equation*}
P'_t(X_i) = 0 \quad \forall t \in \Real.
\end{equation*}
Since $P_t(E_{l,i,\xi}) = \sum_{j=1}^{\alpha}p_{i,j}(t) E_{l,j,\xi}$ for some functions
$p_{i,j}(t)$, we have
\begin{eqnarray*}
&  & P_t'(E_{l,i,\xi} - a_i E_{l,\alpha,\xi}) = P_t'(E_{l,i,\xi}) - a_i P_t'(E_{l,\alpha,\xi}) \\
& = & \sum_{j=1}^{\alpha} p_{i,j}'(t)E_{l,j,\xi} - a_i \sum_{j=1}^{\alpha} p_{\alpha,
  j}'(t)E_{l,j,\xi} = 0.
\end{eqnarray*}
Therefore we have the following system of ordinary differential equations for $p_{i,j}(t)$:
\begin{equation}\label{eqnpijode}
p_{i,j}'(t) - a_i p_{\alpha,j}'(t) = 0. \quad \forall i,j = 1,..., \alpha.
\end{equation}

One easily sees that it has the solution
\begin{equation*}
p_{i,j}(t) = a_i a_j h^2(t) + f_{i,j} - a_i f_{j,\alpha}
\end{equation*}
which finishes our proof. \qedbox

\smallskip

Now we consider the case when the multiplicity of $\mu$ in $\rho$ is one, i.e., $\alpha = 1$. Since
$P_{1,1}(0)$ is symmetric and $\Ad_{\Km}$ equivariant, the off-diagonal terms in the formulas
(\ref{eqnequivarmapcpx}) and (\ref{eqnequivarmapquater}) of equivariant endomorphisms vanish, i.e.,
$P_{1,1}(0) = f_{1,1} I_l$ with some constant $f_{1,1}\in \Real$. From a similar argument as in the
previous case, we have
\begin{thm}\label{thmPtmulone}
If the cohomogeneity one manifold $(M,g)$ has non-negative sectional curvature and the multiplicity
of the class one representation $\mu$ in $\rho$ is one, then for any $\xi = 1, \ldots, n-m$, we
have
\begin{enumerate}
\item $P_t(E_{a,1,\xi}) = f_{1,1} E_{a,1,\xi}$, if $a = 1, \cdots, l-1$;
\item $P_t(E_{l,1,\xi}) = p_{1,1}(t) E_{l,1,\xi}$ and $p_{1,1}(t) = (h^2(t) -1)f_{1,1}^2 + f_{1,1}$.
\end{enumerate}
\end{thm}

\subsection{$\mu$ is of complex or quaternionic type.} We consider the complex case first and the
quaternionic case will follow easily.

Let $l=2p$ and $\beta = 2\alpha$. Recall that $W_i$ is the subspace of $V$ such that $\rho|_{W_i} =
\mu$. Choose an orthonormal basis $\set{e_1, \ldots, e_{2p}}$ of $W_i$ such that $\mu(x)$ has the
form $\left(\begin{smallmatrix} A & -B
\\ B & A \end{smallmatrix}\right)$ with $A, B$ being $p\times p$ matrices and $\mu(H')$ fixes the two
vectors $e_p$ and $e_{2p}$. Under this basis any $\Ad_{\mu(K')}$-equivariant endomorphism has the
block-form $\left(\begin{smallmatrix} aI_p & -b I_p \\ b I_p & a I_p \end{smallmatrix}\right)$as in
(\ref{eqnequivarmapcpx}) with constants $a$, $b \in \Real$. Using the fact that the endomorphism
commutes with the rotation in the $\set{e_p, e_{2p}}$ plane, we may assume that $E_{m,m} \in \Kp$.

Since $\mu$ is of complex type, the $\Ad_{\Km}$-equivariant map $P(0)$ has a block form and the
$(i,j)$-block is given by
\begin{equation*}
\begin{pmatrix}
f_{2i-1,2j-1} I_p & f_{2i,2j-1} I_p \\
f_{2i-1,2j} I_p & f_{2i,2j} I_p
\end{pmatrix},
\end{equation*}
where $f_{a,b}\in \Real$ is constant for $a,b = 1,\ldots, \beta =2\alpha$ and satisfies the
following identities:
\begin{equation}\label{eqnfijidenitycpx}
f_{2i-1,2j-1}= f_{2i,2j}, \quad f_{2i,2j-1} + f_{2i-1,2j} = 0, \quad f_{2i-1,2j-1} = f_{2j-1,2i-1},
\quad f_{2i,2j-1} + f_{2j,2i-1} = 0.
\end{equation}
The last two are due to the fact that $P(0)$ is $Q$-symmetric.

Similar to the case when $\mu$ is of real type, we define $E_{a,i,\xi} = E_{r+(i-1)p + a, m+\xi}$
for $a=1, \ldots, p$, $i=1,\ldots, \beta$ and $\xi= 1,\ldots, n-m$. Then we have
\begin{thm}\label{thmPijcpx}
If the cohomogeneity one manifold $(M,g)$ has non-negative sectional curvature and the class one
representation $\mu$ is of complex type, then for any $i = 1, \ldots, \beta$ and $\xi = 1, \ldots,
n-m$, we have
\begin{enumerate}
\item $P_t(E_{a,i,\xi}) = \sum_{j=1}^{\beta} f_{i,j} E_{a,j,\xi}$, if $a = 1, \cdots, p-1$;
\item $P_t(E_{p,i,\xi}) = \sum_{j=1}^{\beta}p_{i,j}(t) E_{p,j,\xi}$ and $p_{i,j}(t)$ is defined as
\begin{equation}\label{eqnpijcpx}
p_{i,j}(t) = (h^2(t)-1) a_i a_j  + f_{i,j},
\end{equation}
where $a_i = f_{i,\beta}$.
\end{enumerate}
\end{thm}

\smallskip

Now we consider the case when $\mu$ is of quaternionic type and the multiplicity of $\mu$ in $\rho$
is bigger than one. Let $l = 4p$ and $\beta = 4\alpha$.

From the formula (\ref{eqnequivarmapquater}) of equivariant endomorphisms in this case and a
similar argument in the complex case, $P(0) = (f_{a,b}I_p)_{1\leq a, b \leq \beta}$ where
$f_{a,b}\in \Real$ is constant and satisfies the following identities:
\begin{eqnarray}\label{eqnfijidentityquater}
& & f_{4i-3,4j-3} = f_{4i-2,4j-2} = f_{4i-1,4j-1} = f_{4i,4j} = f_{4j,4i} \nonumber\\
& & -f_{4i-2,4j-3} = f_{4i-3,4j-2} = f_{4i,4j-3} = -f_{4i-1,4j} = f_{4j-1,4i} \\
& & -f_{4i-1,4j-3} = -f_{4i,4j-2} = f_{4i-3,4j-1} = f_{4i-2, 4j} = -f_{4j-2,4i} \nonumber \\
& & -f_{4i,4j-3} = f_{4i-1,4j-2} = -f_{4i-2,4j-1} = f_{4i-3,4j} = -f_{4j-3,4i} \nonumber
\end{eqnarray}
for $i,j = 1, \ldots, \alpha$.

We denote $E_{r+(i-1)p+a,m+\xi}$ by $E_{a,i,\xi}$ for $a=1,\ldots, p$, $i=1,\ldots, \beta$ and $\xi
= 1, \ldots, n-m$, and then we have
\begin{thm}\label{thmPijquater}
If the cohomogeneity one manifold $(M,g)$ has non-negative sectional curvature and the class one
representation $\mu$ is of quaternionic type, then for any $i = 1, \ldots, \beta$ and $\xi = 1,
\ldots, n-m$, we have
\begin{enumerate}
\item $P_t(E_{a,i,\xi}) = \sum_{j=1}^{\beta} f_{i,j} E_{a,j,\xi}$, if $a = 1, \cdots, p-1$;
\item $P_t(E_{p,i,\xi}) = \sum_{j=1}^{\beta} p_{i,j}(t) E_{p,j,\xi}$ and $p_{i,j}(t)$ is defined as
\begin{equation}\label{eqnpijcpx}
p_{i,j}(t) = (h^2(t)-1) a_i a_j  + f_{i,j},
\end{equation}
where $a_i = f_{i,\beta}$.
\end{enumerate}
\end{thm}

\section{Proof Of Theorem \ref{thmorth}}\label{secprooforth}

In this section, we will develop some contradiction from the assumption that the manifold $(M,g)$
is non-negatively curved. First we list some lemmas which will be used in the proof of Theorem
\ref{thmorth}.

\begin{lem}\label{lemnontransitive}
The image $\mu(K')\subset SO(l)$ does not act transitively on the sphere
$\sph^{l-1}=SO(l)/SO(l-1)$.
\end{lem}

\pf Recall that $l=\deg \mu \geq k+2$. If $\mu(K')$ could act transitively on $\sph^{l-1}$, then
$K'$ would act transitively on both $\sph^k$ and $\sph^{l-1}$ and $l-1\geq k + 1$. By the
classification of the transitive action on the spheres, we have that $(K',H')$ is either
$(SO(7),SO(6))$ with $\mu(SO(7)) = Spin(7)\subset SO(8)$ or $(SO(9),SO(8))$ with
$\mu(SO(9))=Spin(9)\subset SO(16)$. But in both cases, by the classification in Theorem
\ref{thmclassonecpx}, $\mu$ is not a class one representation for the pair $(K',H')$ which is a
contradiction. \qedbox

\smallskip

The following was already used in \cite{GVWZ} and we state it as a lemma without proof.

\begin{lem}\label{lemanalytic}
Suppose $f(t)$ is a $C^2$ non constant even function on $(-\eps, \eps)$ for some $\eps >0$ with
$f(0) = 0$, then there is no such constant $\gamma \geq 0$ that satisfies the following inequality:
\begin{equation}\label{ineqnder}
\gamma^2(f(t))^2 - (f'(t))^2  \geq 0
\end{equation}
\end{lem}

We will compute the sectional curvatures for some $2$-planes in our examples. The formula in terms
of $P_t$ is well established in \cite{GZposRic} and we quote it for the convenience of the reader.

\begin{thm}[Grove-Ziller]\label{thmformulasec}
If $X, Y \in \frakp$, the sectional curvatures of $M$ at $c(t)$ are determined by
\begin{eqnarray*}
(a)\quad \quad  g(R(X,Y)X,Y) & = & Q(\Am(X,Y),[X,Y]) - \frac{3}{4}Q(P[X,Y]_{\frakp},[X,Y]_{\frakp}) \\
             &   & +Q(\Ap(X,Y), P^{-1}\Ap(X,Y)) - Q(\Ap(X,X), P^{-1}\Ap(Y,Y))\\
             &   & +\frac{1}{4}Q(P'X,Y)^2 -\frac{1}{4}Q(P'X,X)Q(P'Y,Y) \\
(b)\quad \quad  g(R(X,Y)T,Y) & = & -\half Q(P'X,P^{-1}\Ap(Y,Y))+\half Q(P'Y, P^{-1}\Ap(X,Y)) \\
             &   & +\frac{3}{4}Q([X,Y],P'Y) \\
(c)\quad \quad  g(R(X,T)X,T) & = & Q((-\half P'' + \frac{1}{4} P'P^{-1}P')X,X).
\end{eqnarray*}
\end{thm}

\begin{rem}
In the formulas above, $[X,Y]_{\frakp}$ is the $\frakp$-component of $[X,Y]$ and $A_\pm : \frakp
\times \frakp \To \frakg$ are defined as
\begin{equation*}
A_\pm(X,Y) = \half([X, P_t Y]\mp [P_t X, Y]).
\end{equation*}
\end{rem}

\smallskip

\smallskip

From Lemma \ref{lemweylsymmetry}, $h(t)$ is an even function with $h(0)\ne 0$ and $h(L) = 0$, so
$f(t) = h^2(0)-h^2(t) = 1 - h^2(t)$ satisfies the assumptions in Lemma \ref{lemanalytic}. We will
show the inequality (\ref{ineqnder}) holds for some constant $\gamma \geq 0$ using nonnegativity of
the sectional curvature of a carefully chosen $2$-plane. Theorem \ref{thmorth} then follows from
Lemma \ref{lemanalytic}.

\smallskip

\noindent \emph{Proof of Theorem \ref{thmorth}: } In the following argument, only the entries in
the lower right $(n-r)\times (n-r)$-block of $\frakg=\lieso(n)$ are involved, so without loss of
generality we may assume that $r=0$ or equivalently the representation $\rho$ is a sum of $\al$
copies of $\mu$.

\smallskip

\textsc{Case 1:} The representation $\mu$ is of real type. Since $P_0$ defined in
(\ref{eqnPinitial}) is symmetric and positive definite, we can write
\begin{equation*}
\begin{pmatrix}
f_{1,1} & \cdots & f_{1,\al} \\
\vdots  & \ddots & \vdots \\
f_{\al,1} & \cdots & f_{\al,\al}
\end{pmatrix} = A D A^\top,
\end{equation*}
where $A = (A_{i,j})_{\alpha\times \alpha}$ is an orthogonal matrix and $D = \diag (d_1, \ldots,
d_\alpha)$ is a diagonal matrix with positive entries. Define the following vectors in $\frakqsum$:
\begin{eqnarray}
X^u & = & \sum_{i=1}^{l-1}b_iE_{i,u,1} + E_{l,u,2}\label{eqnXu}\\
Y^u & = & \sum_{i=1}^{l-1}b_iE_{i,u,2} + E_{l,u,1},\label{eqnYu}
\end{eqnarray}
where $u = 1, ..., \alpha$ and $\sum_{i=1}^{l-1}b_i^2 = 1$. Further conditions of the $b_i$'s will
be determined later on.

In the matrix $A$, there is a column, say the $i_0$-th column, with $A_{\alpha, i_0}\neq 0$. We
denote $A_{u,i_0}$ by $\Au$, $u = 1,\ldots, \alpha$, and define the following two vectors $X$, $Y$
in $\frakp$:
\begin{equation}
X = \sum_{u=1}^{\alpha} \Au X^u, \quad Y = \sum_{u=1}^{\alpha} \Au Y^u. \label{eqnXY}
\end{equation}

From the definitions of $X^u$ and $Y^v$, it is easy to see that $[X^u, Y^u] = 0$, and if $u \neq
v$, then
\[
[X^u, Y^v] = \sum_{i=1}^{l-1} b_i (E_{vl, i+(u-1)l} + E_{i+(v-1)l,ul}).\]

Hence
\begin{eqnarray*}
[X,Y] & = & \sum_{u,v =1}^{\alpha}\Au \Av [X^u, Y^v]= \sum_{i=1}^{l-1}\sum_{u\neq v}\Au \Av
b_i(E_{vl, i+(u-1)l} +
            E_{i+(v-1)l,ul})\\
      & = & 0.
\end{eqnarray*}

The fact that $X$ commutes with $Y$ makes the computation of the sectional curvature of the
$2-$plane spanned by $X^*$ and $Y^*$ easier since the first two terms in the curvature formula
$(a)$ in Theorem \ref{thmformulasec} vanish. The other four terms are computed in Proposition
\ref{propComputationReal} below. If we plug in the result of each term into the formula of the
sectional curvature, we have
\begin{eqnarray*}
||X^*\wedge Y^*||^2 \sec(X^*,Y^*)_{c(t)} & = & (h^2(t)-1)^2 Q(X_0, P_t^{-1}(X_0)) - (d_{i_0}
A_\alpha)^4  h^2(t)(h'(t))^2\\
                     & - & Q(\Ap (X,X), P_t^{-1}\Ap (X,X)).
\end{eqnarray*}
Here $X_0 \in \frakp$ is specified in Proposition \ref{propComputationReal} and is orthogonal to
$\frakkm$ with respect to $Q$.

Since $P_t$ is positive definite as well as $P_t^{-1}$, we have $Q(\Ap (X,X), P_t^{-1}\Ap
(X,X))\geq 0$. Therefore $\sec(X^*,Y^*) \geq 0$ implies that
\begin{equation*}
(h^2(t)-1)^2 Q(X_0, P_t^{-1}(X_0)) -\left(d_{i_0} A_\alpha \right)^4 h^2(t)(h'(t))^2 \geq 0.
\end{equation*}
The existence of the constant $\gamma \geq 0$ follows from the facts that $A_\alpha \neq 0$ and
$Q(X_0,P_t^{-1}(X_0))$ is bounded from above near $t=0$.

\smallskip

\textsc{Case 2:} The representation $\mu$ is not of real type and the multiplicity of $\mu$ in
$\rho$ equals to $1$. In this case $P_0$ is a scalar matrix and $P_t$ is a diagonal matrix. It is
easy to see that the proof in the previous case works if we choose $X=X^1$ and $Y = Y^1$ in
(\ref{eqnXu}) and (\ref{eqnYu}).

\smallskip

\textsc{Case 3:} The representation $\mu$ is not of real type and the multiplicity of $\mu$ in
$\rho$ is bigger than $1$. Let $p=\half l$ and $\beta = 2\alpha$ if $\mu$ is of complex type and
let $p=\quater l$ and $\beta = 4\alpha$ if it is of quaternionic type. In both cases, we have
$p\geq k+2$. In each case we define the vector $E_{a,i,\xi}$ for $a=1,\ldots, p$, $i =
1,\ldots,\beta$ and $\xi = 1, \ldots, n-m$. Then similarly we can define vectors $X^u$ and $Y^u$
for $u=1,\ldots, \beta$ and use them to define the vectors $X$ and $Y$ as in \textsc{Case 1}. The
formulas of $P_t(X)$ and $P_t(Y)$ are obtained by Theorem \ref{thmPijcpx} and Theorem
\ref{thmPijquater} respectively and the rest of the proof will follow \textsc{Case 1}. Note that
the number of the constants $b_i$'s in $X$ or $Y$ is equal to $p -1$. The existence of the vector
$X_0$ which is orthogonal to $\frakkm$ follows from the fact that $p \geq k+2$. \qedbox

\smallskip

In the case where $\mu$ is real type, the non-vanishing terms in the sectional curvature of the
$2$-plane spanned by $X$ and $Y$ are computed in the following

\begin{prop}\label{propComputationReal}
For the vectors $X$ and $Y$ defined in (\ref{eqnXY}), by choosing a proper value for the constant
$b_i$, we have
\begin{enumerate}
\item There exists some $X_0 \in \frakp$ which is orthogonal to $\frakkm$ with respect to $Q$ such
that $\Ap (X,Y) = (h^2(t)-1)X_0$;
\item $\Ap (X,X) = \Ap (Y,Y)$;
\item $Q(P_t'(X),Y) = 0$;
\item $-\dfrac{1}{4}Q(P_t'(X),X)Q(P_t'(Y),Y) = -\left(d_{i_0} A_\alpha \right)^4
h^2(t)(h'(t))^2$.
\end{enumerate}
\end{prop}

\pf First we compute the endomorphism $P_t$ on $X$ and $Y$. From the defining equations
(\ref{eqnXu}) of $X^u$ and (\ref{eqnYu}) of $Y^u$, we have
\begin{equation}\label{eqnPXu}
P_t(X^u) = \sum_{i=1}^{l-1}\sum_{r=1}^{\alpha}b_i f_{u,r}E_{i,r,1} +
\sum_{r=1}^{\alpha}p_{u,r}(t)E_{l,r,2}
\end{equation}
and
\begin{equation}\label{eqnPYv}
P_t(Y^v) = \sum_{j=1}^{l-1}\sum_{s=1}^{\alpha}b_j f_{v,s}E_{j,s,2} +
\sum_{s=1}^{\alpha}p_{v,s}(t)E_{l,s,1},
\end{equation}
then
\begin{eqnarray*}
[X^u, P_t(Y^v)] & = & [\sum_{i=1}^{l-1}b_iE_{i+(u-1)l,\alfo} + E_{ul,\alft}, \quad
                    \sum_{j=1}^{l-1}\sum_{s=1}^{\alpha}b_j f_{v,s}E_{j+(s-1)l,\alft} +
                    \sum_{s=1}^{\alpha}p_{v,s}(t)E_{sl,\alfo}]\\
              & = & \left(\sum_{i=1}^{l-1} b_i^2f_{v,u} - p_{v,u}(t)\right)E_{\alft,\alfo} +
                    \sum_{i=1}^{l-1}\sum_{s=1}^{\alpha}b_i(f_{v,s}E_{i+(s-1)l,ul}-
                    p_{v,s}(t)E_{i+(u-1)l,sl}),
\end{eqnarray*}
and
\begin{eqnarray*}
[P_t(X^u), Y^v] & = & [\sum_{i=1}^{l-1}\sum_{r=1}^{\alpha}b_i f_{u,r}E_{i+(r-1)l,\alfo} +
                    \sum_{r=1}^{\alpha}p_{u,r}(t)E_{rl,\alft}, \quad \sum_{j=1}^{l-1}b_j E_{j+(v-1)l,\alft} +
                    E_{vl,\alfo}]\\
              & = & \left(\sum_{i=1}^{l-1}b_i^2f_{u,v}-p_{u,v}(t)\right)E_{\alft,\alfo} + \sum_{i=1}^{l-1}
                    \sum_{r=1}^{\alpha}b_i(f_{u,r}E_{vl,i+(r-1)l}-p_{u,r}(t)E_{rl,i+(v-1)l}).
\end{eqnarray*}
Therefore
\begin{equation}\label{eqnApXuYv}
\Ap (X^u,Y^v) = \Half
            \sum_{i=1}^{l-1}\sum_{s=1}^{\alpha}b_i(f_{v,s}E_{i+(s-1)l,ul}+f_{u,s}E_{i+(s-1)l,vl}
            -p_{v,s}(t)E_{i+(u-1)l,sl} - p_{u,s}(t)E_{i+(v-1)l,sl}).\\
\end{equation}
From the above equation, only the terms as $E_{j+(r-1)l,wl}$ have nonzero coefficients in $\Ap
(X,Y)$ and it is denoted by $c_{j,r,w}$. From the formula (\ref{eqnApXuYv}) and the bi-linearity of
$\Ap$, we have
\begin{eqnarray}\label{eqncjrw1}
c_{j,r,w} & = & \Half b_j\left( \sum_{v=1}^{\alpha}(f_{v,r}A_w\Av - p_{v,w}(t)A_rA_v) +
\sum_{u=1}^{\alpha}(f_{u,r}\Au A_w - p_{u,w}A_uA_r)\right)\nonumber\\
          & = & b_j \left(A_w\sum_{v=1}^{\alpha} f_{v,r}A_v - A_r\sum_{v=1}^{\alpha}
                p_{v,w}(t)A_v\right).
\end{eqnarray}
We can compute the terms in (\ref{eqncjrw1}) explicitly as follows,
\begin{equation}\label{eqnc1}
\sum_{v=1}^{\alpha}f_{v,r}A_v = \sum_{v=1}^{\alpha}\sum_{i=1}^{\alpha}A_{v,i}d_i A_{r,i} A_{v,i_0}
= (A^\tau A D A^\tau)_{i_0,r} = d_{i_0} A_{r,i_0} = d_{i_0} A_r
\end{equation}
and
\begin{eqnarray}
\sum_{v=1}^{\alpha}p_{v,w}(t)A_v & = & \sum_{v=1}^{\alpha}(a_v a_w h^2(t) A_v +f_{v,w}A_v -a_v
                                       f_{w,\alpha}A_v)\nonumber \\
                                 & = & d_{i_0}A_w + (a_w h^2(t) -f_{w,\alpha})\sum_{v=1}^{\alpha}a_vA_v \nonumber\\
                                 & = & d_{i_0}A_w + a_w(h^2(t) -1)
                                 \sum_{v=1}^{\alpha}f_{v,\alpha}A_v\nonumber\\
                                 & = & d_{i_0}A_w + d_{i_0}a_w A_\alpha (h^2(t) -
                                 1),\label{eqnc2}
\end{eqnarray}
where the first equality follows Theorem \ref{thmPtrealtype}.

By substituting the new expressions (\ref{eqnc1}) and (\ref{eqnc2}) back in the expression
(\ref{eqncjrw1}) of $c_{j,r,w}$, we have
\begin{equation}\label{eqnc}
c_{j,r,w} = - d_{i_0} A_r A_\alpha a_w b_j (h^2(t) - 1).
\end{equation}

If $r\neq w$, then $E_{j+(r-1)l,wl}$ is orthogonal to $\frakkm$ with respect to $Q$. If $r=w$, by
the formula (\ref{eqnc}) for $c_{j,r,w}$, $c_{j,r,r}$ is a multiple of $b_j$. From the assumption
that the representation $\rho$ is the direct sum of $\mu$ and the embedding in (\ref{eqnrhox}), the
image $\mu_*(v)$(for any $v\in \frakk'$) is block-wise diagonally embedded in $\lieso(m)$. Hence if
for some $r$, the vector $v_r = \sum_{j=1}^{l-1}c_{j,r,r}E_{j+(r-1)l,rl}$ is orthogonal to
$\frakkm$, then all vectors $v_q$'s are orthogonal to $\frakkm$. By Lemma \ref{lemnontransitive} of
the non-transitivity of the action $\mu(K')$ on the sphere $SO(l)/SO(l-1)$ and by choosing the
proper values of $b_i$'s, $v_r$ is orthogonal to $\frakkm$. Therefore $\Ap (X,Y) = (h^2(t) - 1)X_0$
for some $X_0 \in \frakp$ which is orthogonal to $\frakkm$ with respect to $Q$ and $(1)$ is proved.

\smallskip

By (\ref{eqnXu}) and (\ref{eqnPXu}), we have
\begin{eqnarray*}
[X^u, P_t(X^v)] & = & [\sum_{i=1}^{l-1}b_iE_{i+(u-1)l,\alfo} + E_{ul,\alft}, \quad
                    \sum_{j=1}^{l-1}\sum_{s=1}^{\alpha}b_j f_{v,s}E_{j+(s-1)l,\alfo} +
                    \sum_{s=1}^{\alpha}p_{v,s}(t)E_{sl,\alft}]\nonumber \\
              & = & \sum_{i\neq j\mbox{, or }u\neq s}b_i b_jf_{v,s}E_{j+(s-1)l,i+(u-1)l} + \sum_{s\neq
              u}p_{v,s}(t)E_{sl,ul}.
\end{eqnarray*}
By (\ref{eqnYu}) and (\ref{eqnPYv}), we have the same result for $[Y^u, P_t(Y^v)]$, therefore $\Ap
(X,X) = [X,P_t(X)] = [Y,P_t(Y)] = \Ap (Y,Y)$ which proves the formula in $(2)$.

\smallskip

Next we will prove the formulas in $(3)$ and $(4)$ which will finish the proof of the proposition.

By (\ref{eqnPXu}) and the differential equation (\ref{eqnpijode}) for $p_{ij}$ we have
\begin{equation*}
P_t'(X^u) = \sum_{r=1}^{\alpha}p'_{r,u}(t)E_{l,r,2} = 2h(t)h'(t)\sum_{r=1}^{\alpha}a_u a_r
E_{l,r,2},
\end{equation*}
so
\begin{eqnarray*}
Q(P_t'(X^u), Y^v) & = & 2 a_u h(t)h'(t)\sum_{r=1}^{\alpha}a_r \left(\sum_{j=1}^{l-1}b_j
Q(E_{l,r,2},
                  E_{j,v,2}) + Q(E_{l,r,2}, E_{l,v,1})\right)\\
            & = & 0,
\end{eqnarray*}
and then $Q(P_t'(X),Y) = 0$.

By taking the inner product with $X^v$ instead of $Y^v$, we have
\begin{eqnarray*}
Q(P_t'(X^u), X^v) & = & 2 a_u h(t)h'(t)\sum_{r=1}^{\alpha}a_r \left(\sum_{j=1}^{l-1}b_j
Q(E_{l,r,2},
                      E_{j,v,1}) + Q(E_{l,r,2}, E_{l,v,2})\right)\\
                & = & 2a_u a_v h(t)h'(t).
\end{eqnarray*}
Therefore
\begin{eqnarray*}
Q(P_t'(X),X) & = & 2\left(\sum_{u,v =1}^{\alpha}A_uA_v a_u a_v \right)h(t)h'(t) = 2\left(\sum_{u
                  =1}^{\alpha}A_u a_u \right)^2 h(t)h'(t)\\
           & = & 2 \left( d_{i_0} A_\alpha \right)^2 h(t)h'(t),
\end{eqnarray*}
where the last equality follows either from (\ref{eqnc1}) or (\ref{eqnc2}). Similarly for
$Q(P_t'(Y),Y)$ and hence $(4)$ follows, which finishes the proof. \qedbox

\smallskip

\begin{rem}
In the introduction, we pointed out that there are unknown but interesting cases when $m\geq k+2$
and $n=m+1$. The minimal dimension of these manifolds is $15$ when $k=2$, $m=5$ and $n=6$. This
manifold has cohomology ring different from the two $15$ dimensional symmetric spaces $\sph^{15}$
and $SO(8)/(SO(5)\times SO(3))$. The geometry of these examples will be studied in another paper.
\end{rem}

\smallskip

\section{Cohomogeneity One Manifolds for $G=U(n)$ And $Sp(n)$}\label{secunitarysymp}

In this section, we will generalize our examples to the cases where $G = U(n)$ and $Sp(n)$. First
let us state the theorem in each case.

\begin{thm}\label{thmexamplescomplex}
Let $K'/H'=\sph^k$ with $k\geq 2$ and $\rho : K'\To U(m)$ be a faithful representation. Suppose
$\rho$ contains a class one representation $\mu: K' \To U(l)$ of the pair $(K',H')$ with $2l\geq
k+2$ and the multiplicity of $\mu$ in $\rho$ is $1$. For any integer $n\geq m+2$, set $G=U(n)$ and
\begin{eqnarray}
\Km & = & \rho(K')\times U(n-m) \subset U(m)\times U(n-m) \subset U(n) \nonumber\\
\Kp & = & \rho(H')\times U(n-m+1)\subset U(m-1)\times U(n-m+1) \subset U(n) \label{groupdiagramcomplex}\\
H & = & \rho(H')\times U(n-m) \subset U(n),\nonumber
\end{eqnarray}
then the cohomogeneity one manifold $M$ defined by the groups $H\subset \set{\Km,\Kp}\subset G$
does not admit a $G$ invariant metric with non-negative sectional curvature.
\end{thm}

\begin{rem}
Proposition \ref{propdimcomparecpx} lists the complex class one representations which have
dimension smaller than $\half (k+2)$. It shows that if $\mul(\mu,\rho) = 1$, then only the defining
representations of $SU(l)$, $U(l)$, $Sp(l)$ and $Sp(l)\times U(1)$ are excluded by the above
theorem. The cohomogeneity one manifolds from these representations are equivariantly diffeomorphic
to the homogeneous spaces $U(n+1)/\Phi(K')\times U(n-l+1)$ where $\Phi : K' \To U(l)$ is the
defining representation, so they carry non-negatively curved metrics.
\end{rem}

Similar to Lemma \ref{lemnontransitive} in the orthogonal case, we have the following lemma:
\begin{lem}\label{lemnontransitivecpx}
Assume that $K'$, $H'$ and $\mu$ as in Theorem \ref{thmexamplescomplex} with $2l\geq k+2$, then
$\mu(K')$ does not act transitively on $\Cp^{l-1} = U(l)/(U(l-1)\times U(1))$.
\end{lem}

\pf From the classification of the transitive actions on complex projective spaces, see
\cite{Bessegeodesic}, p.195, we only need to check the pair $(SU(2),U(1))$ for $\Cp^1$,
$(U(n),U(n-1))$ for $\Cp^{n-1}$ and $(Sp(n),Sp(n-1))$ for $\Cp^{2n-1}$. In the first case, the
subgroup $U(1) \subset SU(2)$ does not fix any vector in $\Cpx^2$. In the last two cases, we have
$2l=k+1$ which contradicts the assumption on $l$. Therefore the action of $\mu(K')$ on $\Cp^{l-1}$
is non-transitive. \qedbox

\smallskip

\smallskip

\noindent \emph{Sketch of the Proof of Theorem \ref{thmexamplescomplex}: } Suppose $\rho$ has the
decomposition $\tau \oplus \mu $ and $\mu$ is not equivalent to any of subrepresentations in
$\tau$. Let $c(t)$ be the normal geodesic connecting the two singular orbits $\Bpm = G/\Kpm$ with
$c(0)=\ptm \in \Bm$ and $c(L)=\ptp \in \Bp$. Similar to the orthogonal case, the Weyl group is
isomorphic to $\Int_2 \times \Int_2$ and the generators have the following representatives:
\begin{equation*}
\wep =
\begin{pmatrix}
 I_{m-1} &     &           \\
         & -1  &           \\
         &     & I_{n-m}
\end{pmatrix},
\end{equation*}
and
\begin{equation*}
\wem =
\begin{pmatrix}
A_1 &     & & & \\
    &     & A_2 & & \\
    &     &     & \eps & \\
    &     &     &      & I_{n-m}
\end{pmatrix} \mbox{ or }
\begin{pmatrix}
A_1 &     & \\
    & \eps I_{l} & \\
    &     &     I_{n-m}
\end{pmatrix},
\end{equation*}
where $A_1 \in U(r)$, $A_2 \in U(l-1)$ and $\eps = \pm 1$.

In addition to the matrices $\set{E_{i,j}}_{1\leq i\neq j \leq n}$, let $F_{i,j}$ be the symmetric
matrix with $\qi = \sqrt{-1}$ in the $i,j-$ and $j,i-$entries if $i\neq j$ and $\sqrt{2}\qi$ in the
$i,i-$entry if $i=j$. Then $\set{E_{i,j}}$ and $\set{F_{i,j}}$ form an orthonormal basis of the Lie
algebra $\lieu(n)$ of $U(n)$ with $Q = -\half  \Re \Tr$, where $\Re$ takes the real part of a
complex number. Without loss of generality, we may assume that $r=0$, i.e., $\rho=\mu$ and then
$m=l$. Let $\frakp$ be the orthogonal complement of the Lie algebra $\frakh$ of $H$ in the Lie
algebra $\frakg=\lieu(n)$ of $G$ and $\frakq$ be the subspace of $\frakp$ spanned by the vectors
$\set{E_{a,i}}$ and $\set{F_{a,i}}$ for $a=1,\ldots, m$ and $i=m+1,\ldots, n$.

Let $h^2(t) = g(E_{m,m+1}^*,E_{m,m+1}^*)_{c(t)}$ and then we may assume that $h(0) = 1$. From
Schur's Lemma and Wilking's Rigidity Theorem, we have

\begin{prop}\label{propPtCpx}
Suppose that the metric $g$ is non-negatively curved, then we have
\begin{enumerate}
\item $P(E_{a,i}) =E_{a,i} \mbox{ and } P(F_{a,i}) = F_{a,i}$, if $a=1,\cdots, m-1$;
\item $P(E_{m,i}) =h^2(t) E_{m,i} \mbox{ and } P(F_{m,i}) = h^2(t) F_{m,i}$,
\end{enumerate}
where $i=m+1, \ldots, n$.
\end{prop}

From the collapsing of the Killing vector field $E^*_{m,m+1}$ at $\ptp$ and Weyl symmetry at
$\ptm$, $h(t)$ is an even function and $h(L)=0$. Let
\begin{eqnarray*}
X & = & \sum_{i=1}^{m-1} b_i E_{i,m+1} + E_{m,m+2} + \sum_{i=1}^{m-1} c_i F_{i,m+1} + F_{m,m+2} \\
Y & = & \sum_{j=1}^{m-1} b_j E_{j,m+2} + E_{m,m+1} + \sum_{j=1}^{m-1} c_j F_{j,m+2} + F_{m,m+1},
\end{eqnarray*}
where $\sum_{i=1}^{m-1} (b_i^2+c_i^2) = 2$.

A computation shows that $[X,Y] = 0$ and results in the following claim:
\begin{claim}
For properly chosen values of the constants $b_i$ and $c_i$, we have
\begin{enumerate}
\item There exists some $X_0 \in \frakp$ which is orthogonal to $\frakkm$, with respect to $Q$ such
that $\Ap(X,Y) = (h^2(t)-1)X_0$;
\item $\Ap(X,X) = \Ap(Y,Y)$;
\item $Q(P'_t(X),Y) = 0$;
\item $-\dfrac{1}{4}Q(P'_t(X),X)Q(P'_t(Y),Y) = -4 (h(t)h'(t))^2$.
\end{enumerate}
\end{claim}
The existence of $X_0$ follows from the non-transitivity of the $\mu(K')$ action on $\Cp^{l-1} =
U(l)/(U(l-1)\times U(1))$ proved in Lemma \ref{lemnontransitivecpx}. The same argument as in the
orthogonal case shows that the non-negativity of the sectional curvature of the $2-$plane spanned
by $X^*$ and $Y^*$ gives the desired contradiction. This completes the proof in the unitary case.
\qedbox

\smallskip

Finally we discuss the case where $G=Sp(n)$.

\begin{thm}\label{thmexamplesquater}
Let $K'/H'=\sph^k$ with $k\geq 2$ and $\rho : K'\To Sp(m)$ be a faithful representation. Suppose
$\rho$ contains a class one representation $\mu : K' \To Sp(l)$ of the pair $(K',H')$ with $4 l\geq
k+2$ and the multiplicity of $\mu$ in $\rho$ is $1$. For any integer $n\geq m+2$, set $G=Sp(n)$ and
\begin{eqnarray}
\Km & = & \rho(K')\times Sp(n-m) \subset Sp(m)\times Sp(n-m) \subset Sp(n) \nonumber\\
\Kp & = & \rho(H')\times Sp(n-m+1)\subset Sp(m-1)\times Sp(n-m+1) \subset Sp(n) \label{groupdiagramqua}\\
H & = & \rho(H')\times Sp(n-m) \subset Sp(n),\nonumber
\end{eqnarray}
then the cohomogeneity one manifold $M$ defined by the groups $H\subset \set{\Km,\Kp}\subset G$
does not admit a $G$ invariant metric with non-negative sectional curvature.
\end{thm}

\begin{rem}
Proposition \ref{propdimcomparequa} lists the quaternionic class one representations which have
dimension smaller than $\frac{1}{4}(k+2)$. It show that only the standard representation of $Sp(l)$
for the pair $(Sp(l),Sp(l-1))$ is excluded by the above theorem. The cohomogeneity one manifold
from this representation has non-negatively curved metric since it is equivariantly diffeomorphic
to the homogeneous space $Sp(n+1)/Sp(l)\times Sp(n-l+1)$.
\end{rem}

We have the following result on quaternionic projective spaces which is analogues to Lemma
\ref{lemnontransitive} and Lemma \ref{lemnontransitivecpx}:
\begin{lem}\label{lemnontransitivequater}
Assume that $K'$, $H'$ and $\mu$ are as in Theorem \ref{thmexamplesquater} with $4l\geq k+2$, then
$\mu(K')$ does not act transitively on $\Hp^{l-1} = Sp(l)/(Sp(l-1)\times Sp(1))$.
\end{lem}

\pf From the classification of the transitive actions on $\Hp^{l-1}$, see \cite{Bessegeodesic},
p.195, we have $\mu(K') = Sp(l)$ and then $H'=Sp(l-1)$, $\mu$ is the standard representation of $K'
= Sp(l)$. However in this case, $k = 4l -1$ which contradicts the assumption $4l \geq k+2$. \qedbox

\smallskip

\noindent \emph{Sketch of the Proof of Theorem \ref{thmexamplesquater}: } The proof follows as in
the complex case  where $G= U(n)$ step by step. Suppose $\set{1,\qi, \qj, \qk}$ is the basis of
$\Qua$ over the reals such that $\qi^2 = \qj^2 = \qk^2 = -1$ and $\qk = \qi \qj$. Let $G_{i,j}$
denote the symmetric matrix with $1$ in the $i,j-$ and $j,i-$entries if $i \ne j$, and $\sqrt{2}$
in the $i,i-$entries. Thus $\set{E_{i,j}, \qi G_{i,j}, \qj G_{i,j}, \qk G_{i,j}}$ forms an
orthonormal basis of the Lie algebra $\liesp(n)$ of $Sp(n)$ with $Q= - \half \Re \Tr$. Without loss
of generality, we may assume that $r = 0$, i.e., $m = l$.

Let $h^2(t) = g(E^*_{m,m+1},E^*_{m,m+1})_{c(t)}$ and we may assume that $h(0) = 1$. For the
endomorphism $P_t$, one proves the following proposition which is similar to the orthogonal and
complex cases.
\begin{prop}\label{propPtQuater}
Suppose that the metric $g$ is non-negatively curved, then we have
\begin{enumerate}
\item $P(E_{a,i}) = E_{a,i}$ and $P(\theta G_{a,i}) = \theta G_{a,i}$, if $a= 1,\cdots, m-1$;
\item $P(E_{m,i}) = h^2(t) E_{m,i}$ and $P(\theta G_{m,i}) = h^2(t) \theta G_{m,i}$,
\end{enumerate}
where $i = m+1, \ldots, n$ and $\theta$ can be one of $\qi$, $\qj$ and $\qk$.
\end{prop}
Furthermore, $h(t)$ is an even function with $h(L) = 0$. To get the desired contradiction we choose
the vectors $X$ and $Y$ as follows:
\begin{eqnarray*}
X & = & \sum_{i=1}^{m-1} b_i E_{i,m+1} + E_{m, m+2} + \sum_{i=1}^{m-1} c_i G_{i,m+1} +
         c G_{m, m+2}, \\
Y & = & \sum_{j=1}^{m-1} b_j E_{j,m+2} + E_{m, m+1} + \sum_{j=1}^{m-1} c_j G_{j,m+2} +
         c G_{m, m+1},
\end{eqnarray*}
where $c=\qi + \qj +\qk$, $b_i$'s are real numbers and $c_i$'s are pure quaternionic numbers, i.e.
the real part is zero. These constants satisfy the equation $1 - c^2 - \sumi(b_i^2 -c_i^2) = 0$.

A computation shows that $[X,Y] = 0$ and the following claim:
\begin{claim}
For properly chosen values of the constants $b_i$ and $c_i$, we have
\begin{enumerate}
\item There exists some $X_0 \in \frakp$ which orthogonal to $\frakkm$, with respect to $Q$ such
that $\Ap(X,Y) = (h^2(t)-h^2(0))X_0$;
\item $\Ap(X,X) = \Ap(Y,Y)$;
\item $Q(P'_t(X),Y) = 0$;
\item $-\dfrac{1}{4}Q(P'_t(X),X)Q(P'_t(Y),Y) = -16 (h(t)h'(t))^2$.
\end{enumerate}
\end{claim}
The existence of $X_0$ follows from the non-transitivity of the $\mu(K')$ action on $\Hp^{l-1}$
proved in Lemma \ref{lemnontransitivequater}. The contradiction now follows as in the orthogonal
case. \qedbox

\medskip

\medskip

\appendix

\section{Class One Representations of Sphere Group Pairs}\label{appclassone}

In this appendix, we shall classify all class one representations of spherical pairs. Except for
the last section, all representations are considered over the complex numbers.

First we reduce our classification to the case where the group action on the sphere is almost
effective. Suppose $K$, $H$ are compact Lie groups, $K$ is connected and $K/H=\sph^{n-1}$ with $n
\geq 2$. If the $K$ action is not almost effective, then let $C$ be the ineffective kernel, i.e.
the maximal normal subgroup of $K$ contained in $H$. Thus we can write $K=C\times K_1$ and $H=
C\times H_1$. Suppose $\mu \otimes \tau$ is a class one representation of the pair $(K,H)$, i.e.,
$\mu$ and $\tau$ are irreducible representations of $C$ and $K_1$ respectively and $\Res(\mu
\otimes \tau)$ fixes a non-zero vector, i.e., the trivial representation of $C\times H_1$ appears
in the decomposition of $\mu \otimes \tau$. Therefore $\mu$ is the trivial representation of $C$
and $\tau$ is a class one representation of the pair $(K_1, H_1)$. Therefore in the rest, we only
consider almost effective action on the spheres.

The classification of transitive and effective action on spheres by connected compact Lie groups is
well known, see for example Page $179$ in \cite{BesseEinstein} or Page $195$ in
\cite{Bessegeodesic}. Using representation theory of compact Lie groups, it is easy to extend the
classification to the almost effective case:
\begin{itemize}
  \item $SO(n)/SO(n-1) = \sph^{n-1}$($n\geq 3$),
  \item $SU(n)/SU(n-1) = \sph^{2n-1}$($n\geq 2$),
  \item $U(n)/U(n-1)_m = \sph^{2n-1}$($n\geq 1 \mbox{ and } m\in \Int$),
  \item $Sp(n)/Sp(n-1) = \sph^{4n-1}$($n\geq 1$),
  \item $Sp(n)\times Sp(1) / (Sp(n-1)\times Sp(1)) = \sph^{4n-1}$($n \geq 1$),
  \item $Sp(n)\times U(1) / (Sp(n-1)\times U(1))_m = \sph^{4n-1}$($n\geq 1$ and $m \in \Int$, $m \ne 0$),
  \item $G_2/SU(3) = \sph^6, \quad Spin(7)/G_2 = \sph^7, \quad Spin(9)/Spin(7) = \sph^{15}$.
\end{itemize}

\smallskip

The group $U(n)$ can act on the sphere $\sph^{2n-1}$ in different ways. For each $m \in \Int$, $A
\in U(n)$ can act on $\Cpx^n$ via $z \mapsto (\det A)^m A.z$. It induces a transitive action on
$\sph^{2n-1} \subset \Cpx^n$ and the isotropy subgroup at $z=(1, 0, \cdots, 0)^\top$ is
\begin{equation*}
U(n-1)_m = \set{\diag(a, B) \in U(n) | a^{m+1} = (\det B)^{-m}, a\in U(1) \mbox{ and } B \in
U(n-1)}.
\end{equation*}

If $n=1$, then $U(1)$ acts on the circle $\sph^1$ and the isotropy subgroup is $\Int_{m+1}(m\ne
-1)$ and it gives us all almost effective actions on the circle. In Table \ref{tableclassonecpx},
we only list the case when the $U(1)$ action on $\sph^1$ is effective and the corresponding results
in the almost effective case easily follows.

Similarly the group $Sp(n) \times U(1)$ has different transitive actions on $\sph^{4n-1}$. For each
$m \in \Int$, $(A,z) \in Sp(n) \times U(1)$ acts on $\Qua^n$ via $q \mapsto A.q z^m$. The isotropy
subgroup at $q=(0,\cdots, 0, 1)^\top$ is
\begin{equation*}
Sp(n-1)\times U(1)_m = \set{(\diag(A,z^{-m}),z) \in Sp(n) \times U(1) | z \in U(1) \mbox{ and }A
\in Sp(n-1)}.
\end{equation*}

\smallskip

For each spherical pair $(K,H)$ with $K/H = \sph^{n-1}$, the defining representation $\Phi : K \To
SO(n)$ is of class one. If the pair is $(SO(n), SO(n-1))$, then the class one representations are
well known. They consist of the irreducible representations on the space of homogeneous harmonic
polynomials. In fact the class one representations of the spherical pairs are closely related to
these representations as stated in the following theorem.
\begin{thm}\label{thmharmonic}
The representation $\mu$ is a class one representation of the spherical pair $(K,H)$ if and only if
$\mu$ is in the decomposition of $\Res^{SO(n)}_{K} \rho$, where $\rho$ is a class one
representation of the pair $(SO(n), SO(n-1))$ and $K$ is viewed as a subgroup of $SO(n)$ via the
representation $\Phi$.
\end{thm}

The proof of the above theorem is given in Section \ref{secclassoneharmonic}.

\smallskip

For a compact Lie group, the irreducible representations are highest weight representations and
each highest weight is a linear combination of the fundamental weights with nonnegative integers as
coefficients. We list the fundamental weights for classical groups as follows.
\[
\begin{array}{rcl}
SO(n) & : & \fw_1 = e_1, \cdots, \fw_{k-1} = e_1 + e_2 + \cdots + e_{k-1}, \fw_k = \half(e_1 + \cdots + e_k), \mbox{ n = 2k+1,} \\
      &   & \fw_1 = e_1, \cdots, \fw_{k-1} = \half(e_1 + \cdots + e_{k-1} -e_k), \fw_{k} = \half(e_1 + \cdots + e_k), \mbox{ n = 2k},\\
SU(n) & : & \fw_1 = e_1, \cdots, \fw_{n-1} = e_1 + \cdots + e_{n-1}, \mbox{ with } e_1 + \cdots + e_n = 0,\\
Sp(n) & : & \fw_1 = e_1, \cdots, \fw_n = e_1 + \cdots + e_n.
\end{array}
\]
The exceptional Lie group $G_2$ has two fundamental weights: $\fw_1$ which is the highest weight of
the $7$ dimensional representation and $\fw_2$ which is $14$ dimensional.

The group $U(n)$ has a finite cover $SU(n)\times U(1)$ and hence its irreducible representations
can be written as $\mu \otimes \phi^k$ where $\mu$ is an irreducible representation of $SU(n)$ and
$\phi^k$($k \in \Int$) is the $1$ dimensional representation of $U(1)$: $v \mapsto z^k v$ for any
$z\in U(1)$. Hence an irreducible representation $\rho$ of $U(n)$ with highest weight $a_1 e_1 +
\ldots + a_n e_n$( $a_1 \geq \ldots \geq a_n$) is the tensor product of an irreducible
representation $\mu$ of $SU(n)$ and $\phi^k$ of $U(1)$ where $\mu$ has highest weight $(a_1
-a_n)e_1 + \ldots + (a_{n-1} - a_n)e_{n-1}$ and $k = -(a_1 + \ldots + a_n)$. Note that the standard
representation of $U(n)$ on $\Cpx^n$, i.e., matrix multiplication from the left, has the highest
weight $-e_n$.

\smallskip

In order to determine the class one representations for each pair $(K,H)$, we use the branching
rules, i.e. the rules that show how an irreducible representation of $K$ decomposes under the
restriction functor $\Res^K_H$. If the trivial representation of $H$ appears in the decomposition,
then this representation is of class one. We have the following classification result:
\begin{thm}\label{thmclassonecpx}
For each almost effective spherical pair $(K,H)$, Table \ref{tableclassonecpx} gives the
classification of all complex irreducible class one representations $\rho$.

The multiplicity of the trivial representation of $H$ in $\Res^K_H(\rho)$ is equal to $1$ except
for the pair $(Sp(n),Sp(n-1))$ where it is $a+1$ if $\rho$ has the highest weight $a\fw_1 +
b\fw_2$.

In the table the numbers $a$ and $b$ are non-negative integers and $k$ is an integer. In the pair
$(Sp(n)\times U(1), Sp(n-1)\times U(1)_m)$, $a$, $b$ and $k$ satisfy further restriction, denoted
by (S): $a+b \geq 1$, $\abs{m}$ divides $k$, $a$ and $\dfrac{k}{\abs{m}}$ have the same parity and
$\dfrac{\abs{k}}{\abs{m}} \leq a$. For other pairs, the restrictions are specified in the table.

\begin{table}[!h]
\begin{center}
\begin{tabular}{|c|c|c|c|c|c|}
  \hline
  $K$  &  $H$   &  $\rho$                           &  & $n$\\
  \hline
  \hline
  $SO(n)$ &  $SO(n-1)$  & $a \fw_1$ & $a\geq 1$ & $n\geq 3$ \\
  \hline
  $SU(n)$ &  $SU(n-1)$  & $a \fw_1 + b\fw_{n-1}$ & $a+b \geq 1$ & $n\geq 3$ \\
  \hline
  $U(n) $ &  $U(n-1)_m$  & $ae_1 -b e_n + m(a-b)(e_1 + \cdots + e_n)$ & $a+b\geq 1$ & $n \geq 2$ \\
  \hline
  $Sp(n)$ & $Sp(n-1)$ & $a\fw_1 + b\fw_2$ & $a+b \geq 1$ & $n \geq 1$ \\
  \hline
  \hline
  $Sp(n)\times Sp(1)$ &  $Sp(n-1)\times Sp(1)$  & $(a \fw_1 + b \fw_2)\otimes a \fw_1$ &
  $a+b \geq 1$ & $n \geq 1$\\
  \hline
  $Sp(n)\times U(1)$ &  $Sp(n-1)\times U(1)_m$  & $(a \fw_1 + b\fw_2) \otimes \phi^k$ &
  (S) & $n\geq 1$\\
  \hline
  \hline
  $G_2$ &  $SU(3)$  & $a \fw_1$ & $a\geq 1$ & \\
  \hline
  $Spin(7) $ &  $G_2$  & $a\fw_3 $ & $a\geq 1$ & \\
  \hline
  $Spin(9) $ &  $Spin(7)$  & $a\fw_1 + b\fw_4$ & $a+b \geq 1$ & \\
  \hline
  \hline
  $U(1)$ & $\set{1}$  & $\phi^k$ & $k\ne 0$ &  \\
  \hline
\end{tabular}
\vskip 0.1cm \caption{Complex class one representations of spherical pairs}
\label{tableclassonecpx}
\end{center}
\end{table}
\end{thm}

In Table \ref{tableclassonecpx}, if $n=1$ for the $Sp(n)$ factor, then there does not exist
$b\fw_2$ in $\rho$, i.e., $b=0$.

\smallskip

The proof of Theorem \ref{thmclassonecpx} is divided into four parts. The results of the first four
group pairs in Table \ref{tableclassonecpx} are direct consequences of the classical branching
rules for those pairs. The second part includes the pairs $(Sp(n)\times Sp(1), Sp(n-1)\times
Sp(1))$ and $(Sp(n)\times U(1), Sp(n-1)\times U(1)_m)$. They will be proved in the Section
\ref{secsp} by using a branching rule of J. Lepowsky. The third part is covered in Section
\ref{secexceptional}. We will apply Kostant's Branching Theorem to each pair to obtain the class
one representations. The result for the last pair $(U(1),\set{1})$ is clear.

In the last section, Section \ref{secproperties}, we will study some properties of these
representations, for example, the type, the kernel and the dimension.

\subsection{The Class One Representations and Harmonic Polynomials}\label{secclassoneharmonic}

In this section we prove Theorem \ref{thmharmonic} and provide an explicit characterization of
class one representations for the pairs $(SU(n),SU(n-1))$ and $(U(n),U(n-1)_m)$ using harmonic
polynomials.

\smallskip

\emph{Proof of Theorem \ref{thmharmonic}} : We consider the separable Hilbert space $\Hil$ of all
the square-integrable complex valued functions defined on $\sph^{n-1}$ under the following
inner-product
\[
<f_1,f_2> = \int_{\sph^{n-1}}f_1\bar{f_2} d\sigma, \quad \forall f_1, f_2 \in \Hil,
\]
where $d\sigma$ is an $SO(n)$ bi-invariant measure on $\sph^{n-1}$ with total measure $= 1$.

Using the $SO(n)$ action on $\sph^{n-1}$, $SO(n)$ acts on $\Hil$ by $g\star f(x) = f(g^{-1}x)$ for
any $g \in SO(n)$. This representation is called the \emph{left-regular representation} of $SO(n)$
on $\Hil$ and it is a unitary representation.

Since $SO(n)/SO(n-1) = \sph^{n-1}$ is a homogeneous space, $\Hil$ is also the induced
representation $\Hil = \Ind^{SO(n)}_{SO(n-1)}\Id$, where $\Id$ is the trivial representation of
$SO(n-1)$ on $\Cpx$, see for example, Chapter 9.2 in \cite{Knappbook}. By the Frobenius
reciprocity, for any irreducible representation $\rho$ of $SO(n)$, we have the following
multiplicities equality
\[
[\Ind^{SO(n)}_{SO(n-1)}\Id : \rho ] = [\Res^{SO(n)}_{SO(n-1)}\rho : \Id], \quad \mbox{or}\quad
[\Hil : \rho] = [\Res^{SO(n)}_{SO(n-1)}\rho : \Id],
\]
where $[\rho_1 : \rho_2]$ denotes the multiplicity of the irreducible representation $\rho_2$ in
$\rho_1$. By the classical branching rule for the pair $(SO(n), SO(n-1))$,
$[\Res^{SO(n)}_{SO(n-1)}\rho : \Id] \neq 0$ if and only if $\rho$ has the highest weight $ae_1$
where $a$ is a nonnegative integer. On the other hand, for any class one representation $\rho$,
$[\Res^{SO(n)}_{SO(n-1)}\rho : \Id] =1$, or equivalently $[\Hil : \rho] = 1$.

The representation with highest weight $ae_1$ can be realized as the the representation of the
complex valued homogeneous harmonic polynomials with degree $a$ on $\sph^{n-1}$ which is denoted by
$\mathcal{H}^a$. Furthermore, $\mathcal{H}^a$ is the complexification of the real valued
homogeneous harmonic polynomials with degree $a$. We have the following orthogonal decomposition
\[
\Hil = \bigoplus_{a=0}^{\infty}\mathcal{H}^a.
\]

Since $K/H = \sph^{n-1}$, $K$ acts on the Hilbert space $\Hil$ and $\Hil = \Ind^K_H \Id$, where
$\Id$ is the trivial representation of $H$ on $\Cpx$. As a subgroup of $SO(n)$, $K$ acts
invariantly on $\mathcal{H}^a$. In general $\mathcal{H}^a$ is not an irreducible representation of
$K$ and it is decomposed orthogonally into irreducible summands as follows:
\[
\Res^{SO(n)}_K \mathcal{H}^a = \bigoplus_b \mathcal{H}^{a,b},
\]
and there are only finite many $b$'s for each value of $a$. Hence we have the following orthogonal
decomposition of $\Hil$ into irreducible representations of $K$:
\[
\Hil = \bigoplus_{a,b}\mathcal{H}^{a,b}.
\]
Suppose $\mu$ is an irreducible representation of $K$, then $\mu$ is of class one for the pair
$(K,H)$ if and only if $[\Hil : \mu] \neq 0$ by the Frobenius reciprocity applied to the pair
$(K,H)$. Then from the above decomposition and the uniqueness of this decomposition, $\mu$ is
equivalent to $\mathcal{H}^{a,b}$. \qedbox

\begin{rem}
The existence and uniqueness of the decomposition of $\Hil$ follow from Theorem 9.4 and Corollary
9.6 in \cite{Knappbook}.
\end{rem}

\smallskip

Next using harmonic polynomials, we construct the representation space of the class one
representations $\rho = a\fw_1 + b\fw_{n-1}$ for the pair $(K,H)=(SU(n),SU(n-1))$ and $\rho= a e_1
-b e_n + m(a-b)(e_1 + \cdots + e_n)$ for the pair $(K,H)=(U(n),U(n-1)_m)$. It is already discussed
for the pair $(SU(n),SU(n-1))$ in \cite{Knappbook}. In the following, we consider the pair
$(U(n),U(n-1)_m)$.

Let $\set{z_1, \cdots, z_n}$ be the basis of $\Cpx^n$ and $A \in K$ acts on it by the defining
representation $\Phi_m(A) : (z_1, \cdots, z_n)^\top \mapsto (\det A)^m A.(z_1, \cdots, z_n)^\top $.
If $A\in SU(n)$, then $\Phi_m(A)$ is just matrix multiplication. Similarly $\Phi_m(A)$ maps
$(\bar{z}_1, \cdots, \bar{z}_n)^\top$ to $(\det \bar{A})^{m}\bar{A}.(\bar{z}_1, \cdots,
\bar{z}_n)^\top$. Let $V$ be the space of homogeneous polynomials in $z_1, \cdots, z_n$,
$\bar{z}_1, \cdots, \bar{z}_n$ of degree $a+b$. $A$ acts on $V$ by
\begin{equation*}
\Phi_m(A) : P(
\begin{pmatrix}
z_1 \\
\vdots \\
z_n
\end{pmatrix},
\begin{pmatrix}
\bar{z}_1 \\
\vdots \\
\bar{z}_n
\end{pmatrix}
) \mapsto P((\det A)^{-m}A^{-1}.
\begin{pmatrix}
z_1 \\
\vdots \\
z_n
\end{pmatrix},(\det\bar{A})^{-m}\bar{A}^{-1}.
\begin{pmatrix}
\bar{z}_1 \\
\vdots \\
\bar{z}_n
\end{pmatrix}
).
\end{equation*}

Let $V_{a,b}$ be the subspace of polynomials with degree $a$ in $\bar{z}$ and degree $b$ in $z$.
Clearly $K$ leaves $V_{a,b}$ invariant. The Laplacian operator $\Delta$ is a multiple of
$\sum_{j=1}^{n}\frac{\partial^2}{\partial z_j \partial \bar{z}_j}$ and it commutes with the $K$
action, so the subspace $H_{a,b}$ of harmonic polynomials in $V_{a,b}$ is an invariant subspace.

If $P$ is a monomial of the form
\begin{equation*}
P(z, \bar{z}) = \bar{z}_1^{k_1}\cdots \bar{z}_n^{k_n}z_1^{l_1}\cdots z_n^{l_n} \quad \mbox{with }
\sum_{j=1}^{n}k_j = a \mbox{ and } \sum_{j=1}^{n}l_j =b,
\end{equation*}
then $P$ is a weight vector of weight $\sum_{j=1}^n(k_j - l_j)e_j + m(a-b)(e_1 + \cdots + e_n)$.
Hence $\rho = a e_1 - b e_n + m(a-b)(e_1 + \cdots + e_n)$ is the highest weight of the
representation $H_{a,b}$ with the weight vector $\bar{z}_1^a z_n^b$. From Exercise $17$ of Chapter
IV in \cite{Knappbook}, the dimension of $H_{a,b}$ is
$\frac{a+b+n-1}{n-1}\binom{a+n-2}{a}\binom{b+n-2}{n-2}$ which equals to the dimension of the
irreducible representation of $K$ with the highest weight $\rho$, see, for example, Proposition
\ref{propclassonedim}, so $H_{a,b}$ is the representation space as desired.

\smallskip

\subsection{The Pairs $(Sp(n)\times Sp(1), Sp(n-1)\times Sp(1))$ and $(Sp(n)\times U(1), Sp(n-1)\times
U(1)_m)$}\label{secsp}

For each pair $(K,H)$, there is an intermediate group $L$ and we will apply the branching rule
successively as $\Res^{K}_{H}\rho = \Res^L_H (\Res^K_L \rho)$ for any irreducible representation
$\rho$ of $K$.

The branching rule we will use is the one for the pair $(Sp(n),Sp(n-1)\times Sp(1))$ which is
established by J. Lepowsky in \cite{LepowskyMul}. For each irreducible representation $\mu_1
\otimes \mu_2$ of $Sp(n-1)\times Sp(1)$, let $\mu = b_1 e_1 + \cdots + b_{n-1}e_{n-1} + b_n e_n$ be
the highest weight, then $b_1 \geq \cdots \geq b_{n-1}\geq 0$ and $b_n \geq 0$.
\begin{defn}
Let $l,m \in \Int$, $m \geq 1$ and let $q_1, \cdots, q_m$ be positive integers. Suppose there are
$l$ identical balls, $m$ different boxes and the $i$-th box can contain at most $q_i$ balls. We
define $F_m(l;q_1, \cdots, q_m)$ to be the number of ways of putting the $l$ balls into the $m$
boxes.
\end{defn}

\begin{thm}[Lepowsky]\label{thmbranchingSp}
Let $\rho$, $\mu$ be irreducible representations of $Sp(n), Sp(n-1)\times Sp(1)$ with highest
weights $\rho=a_1e_1 + \cdots a_ne_n$, $\mu =b_1e_1 + \cdots b_{n-1}e_{n-1} + b_ne_n$ respectively.
Let
\begin{eqnarray*}
A_1 & = & a_1 -\max(a_2, b_1), \\
A_i & = & \min(a_i,b_{i-1}) - \max(a_{i+1},b_{i}) \quad \quad (2\leq i \leq n-1), \\
A_n & = & \min(a_n,b_{n-1})
\end{eqnarray*}
Then the multiplicity $[\Res\rho : \mu] = 0$ unless $\sum_{i=1}^{n}(a_i+b_i) \in 2\Int$ and $A_1,
\cdots, A_{n-1}\geq 0$. Under these conditions,
\begin{eqnarray}
[\Res\rho : \mu] & = & F_{n-1}\left(\half(b_n - A_1 + \sum_{i=1}^n A_i); A_2, \cdots, A_n\right) \nonumber \\
                 & - & F_{n-1}\left(\half(-b_n -A_1 + \sum_{i=2}^n A_i)-1; A_2, \cdots, A_n\right).\label{eqnmulsp}
\end{eqnarray}
\end{thm}

\smallskip

\emph{Proof of the 2nd part of Theorem \ref{thmclassonecpx} :} Let $K=Sp(n)\times Sp(1)$ and
$H=Sp(n-1)\times Sp(1)$, then $L=Sp(n-1)\times Sp(1)\times Sp(1)$ lies between $H$ and $K$, and the
$Sp(1)$ factor in $H$ is diagonally embedded in $Sp(1)\times Sp(1) \subset L$. Suppose $\rho =
\rho_1\otimes \rho_2$ is a class one representation of the pair $(K,H)$ and $\rho_1$ has the
highest weight $a_1 e_1 + \cdots + a_n e_n$ and $\rho_2$ has the highest weight $be_1$. Then
$\Res^K_L \rho = (\Res^{Sp(n)}_{Sp(n-1)\times Sp(1)} \rho_1) \otimes \rho_2$. Since the trivial
representation appears in $\Res^K_H \rho$, it appears in $(\Res^L_H\mu) \otimes \rho_2$ for some
irreducible representation $\mu = \mu_1\otimes \mu_2$ in $\Res^{Sp(n)}_{Sp(n-1)\times
Sp(1)}\rho_1$. Under the restriction from $Sp(n-1)\times Sp(1)\times Sp(1)$ to $Sp(n-1)\times
\Delta Sp(1)$, $\mu_1 \otimes \mu_2 \otimes \rho_2$ splits as $\mu_1 \otimes \Res^{Sp(1)\times
Sp(1)}_{\Delta Sp(1)} (\mu_2 \otimes \rho_2)$. Therefore $\mu_1$ is a trivial representation of
$Sp(n-1)$ and $\mu_2$, $\rho_2$ have the same dimension, i.e., $\Res^{Sp(n)}_{Sp(n-1)\times
Sp(1)}\rho_1$ contains the representation $\Id\otimes \rho_2$ which has the highest weight $be_n$.
From Theorem \ref{thmbranchingSp}, we have $A_1 = a_1 -a_2 \geq 0$, $A_2 = -a_3 \geq 0$, $\ldots$,
$A_{n-1}= -a_n\geq 0$, $b + a_1 + a_2$ is an even number and the multiplicity is equal to
\begin{equation*}
p = F_{n-1}\left( \half (b-(a_1 -a_2)); 0, \cdots, 0\right) - F_{n-1}\left(\half (-b - (a_1
-a_2))-1; 0, \cdots, 0 \right).
\end{equation*}
Hence $p \ne 0$ if and only if $b=a_1 -a_2$ and in this case we have $p=1$.

For the second part, we give a proof when $m=1$ and the general case easily follows from this case.
Let $K= Sp(n)\times U(1)$, $H = Sp(n-1)\times U(1)$, $K_1 = Sp(n-1)\times Sp(1)\times U(1)$ and
$K_2 = Sp(n-1)\times U(1) \times U(1)$, then we have the embedding $H \subset K_1 \subset K_2
\subset K$ and $U(1)$ in $H$ lies diagonally in $U(1)\times U(1) \subset K_1$. Suppose $\rho = \mu
\otimes \phi^k$ is a class one representation of the pair $(K,H)$ and $\mu$ has the highest weight
$\mu = a_1 e_1 + \cdots + a_n e_n$. From the first part of the proof, we have $a_3 = \cdots =
a_n=0$ and the multiplicity of the trivial representation in $\Res^K_H \rho$ is equal to the
multiplicity of the trivial representation in $(\Res^{Sp(1)}_{U(1)} (a_1-a_2)e_n) \otimes \phi^k$.
Therefore $k$, $a_1 -a_2$ have the same parity, $\abs{k}\leq a_1 -a_2$ and the multiplicity is $1$.
\qedbox

\smallskip

\subsection{The Pairs $(G_2, SU(3))$, $(Spin(7), G_2)$ and
$(Spin(9),Spin(7))$}\label{secexceptional}

In this section, we will develop the branching rule for each group pair and then prove the
corresponding result in Theorem \ref{thmclassonecpx}. The main tool is Kostant's Branching Theorem.
We quote the statement of the theorem from \cite{Knappbook}.

First we establish the notation we will use. Let $K$ be a connected compact Lie group and let $H$
be a connected closed subgroup. Choose a maximal torus $S \subset H$. The special assumption is
that the centralizer of $S$ in $K$ is a maximal torus $T$ of $K$. Let $\Delta_K$ be the set of the
roots of $(\frakk^\Cpx, \frakt^\Cpx)$, let $\Delta_H$ be the set of roots of $(\frakh^\Cpx,
\fraks^\Cpx)$, and let $W_K$ be the Weyl group of $\Delta_K$. Introduce compatible positive systems
$\Delta_K^+$ and $\Delta_H^+$, let bar denote restriction from the dual $(\frakt^\Cpx)^*$ to the
dual $(\fraks^\Cpx)^*$, and let $\delta_K$ be half the sum of the members in $\Delta_K^+$. The
restrictions to $\fraks^\Cpx$ of the members of $\Delta_K^+$, repeated according to their
multiplicities, are the nonzero positive weights of $\fraks^\Cpx$ in $\frakk^\Cpx$. Deleting from
this set the members of $\Delta_H^+$, each with multiplicity $1$, we obtain the set $\Sigma$ of
positive weights of $\fraks^\Cpx$ in $\frakk^\Cpx/\frakh^\Cpx$, repeated according to
multiplicities. The associated Kostant partition function is defined as follows: $\mathcal{P}(\nu)$
is the number of ways that a member of $(\fraks^\Cpx)^*$ can be written as a sum of members of
$\Sigma$, with the multiple versions of a member of $\Sigma$ being regarded as distinct.

\begin{thm}[Kostant's Branching Theorem]\label{thmKostant}
Let $K$ be a compact connected Lie group, let $H$ be a closed connected subgroup, suppose that the
centralizer in $K$ of a maximal torus $S$ of $H$ is abelian and is therefore a maximal torus $T$ of
$K$, and let other notation be as above. Let $\rho$ be an irreducible representation of $K$ with
the highest weight $\rho \in (\frakt^\Cpx)^*$, and let $\mu$ be an irreducible representation of
$H$ with the highest weight $\mu \in (\fraks^\Cpx)^*$. Then the multiplicity of $\mu$ in the
restriction of $\rho$ to $H$ is given by
\[
 [\rho : \mu] = \sum_{w \in W_K}\eps(w)\mathcal{P}(\overline{w(\rho+\delta_K)-\delta_K}-\mu).
\]
\end{thm}

\begin{rem}
To apply this theorem to our special examples, we will use an equivalent assumption on the Cartan
subalgebras instead of the maximal tori: the centralizer of $\fraks$ in $\frakh$ is a Cartan
subalgebra $\frakt$ of $\frakk$. In our three group pairs, this assumption is verified.
\end{rem}

To apply Kostant Branching Theorem, we need to work out the explicit embeddings of the Lie
algebras. These groups and Lie algebras are well studied, see, for examples,
\cite{GluckWarnerZiller} and \cite{Murakami} and references given there.

We start with the pair $(Spin(9),Spin(7))$. Let $\Cay$ be the set of the Cayley numbers which is
isomorphic to $\Real^8$ as vector spaces. $SO(8)$ acts on $\Cay$ by left multiplication and we have
the \emph{Principle of Triality}:
\begin{prop}\label{propTriality}
For each $\theta_1 \in SO(8)$, there exists $\theta_2$, $\theta_3 \in SO(8)$ such that
\begin{equation}\label{eqnTriality}
\theta_1(x)\theta_2(y) = \theta_3(xy), \quad \mbox{for any }x, y \in \Cay
\end{equation}
Moreover if $\theta'_2$, $\theta'_3$ satisfy the above equation as $\theta_2$, $\theta_3$, then
$(\theta'_2, \theta'_3) = \pm (\theta_2, \theta_3)$.
\end{prop}

The infinitesimal version of the above principle is given as follows:
\begin{prop}\label{propTrialityso}
For any $X \in \lieso(8)$, there exist $Y$, $Z \in \lieso(8)$ such that
\begin{equation}\label{eqnTrialityso}
(Xx)y + x(Yy) = Z(xy) \quad \mbox{ for any } x, y \in \Cay.
\end{equation}
Moreover $Y$, $Z$ is uniquely determined by $X$ with $Y = \lambda(X), Z = \lambda\kappa(X)$ and
$\lambda$, $\kappa$ are two outer automorphisms of $\lieso(8)$ with $\lambda^3 = 1$, $\kappa^2 = 1$
and $\kappa \lambda^2 = \lambda \kappa$.
\end{prop}

We identify $\Cay \oplus \Cay$ with $\Real^{16}$ and then $Spin(9)\subset SO(16)$ acts transitively
on $\sph^{15} \subset \Cay \oplus \Cay$. Let $v_0\in \Real \subset \Cay$ be the unit length vector,
then $Spin(7)$ is the isotropy subgroup at $(v_0, 0)$. Consider the following Hopf fibration:
\[
\begin{array}{ccccc}
\sph^7 &  \To & \sph^{15} & \To & \sph^8 = \Cay \cup \set{\infty}, \\
       &      & (x, y)    & \mapsto & y^{-1}\bar{x}.
\end{array}
\]
The isotropy subgroup of the fiber $\set{(x,0)| x\in \Cay}$ is given as follows:
\begin{equation}\label{eqnspin8}
Spin(8) = \set{(\theta_1,\theta_2)\in SO(8)\times SO(8) | \exists  \theta_3 \in SO(8) \mbox{ such
that }\theta_1(x)\theta_2(y) = \overline{\theta_3(\overline{xy}) }},
\end{equation}
and the embedding of $Spin(7)\subset Spin(8)$ is
\begin{equation}\label{eqnspin7}
Spin(7) = \set{(\theta_1,\theta_2)\in Spin(8)| \theta_1 \in SO(7) \subset SO(8) \mbox{ or
}\theta_1(x)\theta_2(y) = \overline{\theta_2(\overline{xy}) }}.
\end{equation}

The automorphism $\theta \mapsto (x \mapsto \overline{\theta(\bar{x})})$ of $SO(8)$ induces the
automorphism $\kappa$. Let $\theta_1 = \exp(X)$, $\theta_2 = \exp(Y)$ and $\theta_3 = \exp(Z)$.
From Proposition \ref{propTrialityso}, we have $X = \lambda(\kappa Z)$ and $Y = \lambda^2(\kappa
Z)$, so we can write the Lie algebra of $Spin(8)$ as
\begin{equation*}
\Lie(Spin(8)) = \set{(\lambda(X),\lambda^2(X))| X \in \lieso(8)}.
\end{equation*}
From the embedding of $Spin(7)\subset Spin(8)$ in (\ref{eqnspin7}), since $\theta_1 \in SO(7)$ we
have the Lie algebra of $Spin(7)$ as
\begin{equation*}
\Lie(Spin(7)) = \set {(Y, \lambda(Y))| Y \in \lieso(7) }\subset \Lie(Spin(8)).
\end{equation*}
Therefore if we identify $\Lie(Spin(8))$ with $\lieso(8)$ by its first component, then we have
\begin{equation}\label{eqnliespin7}
\Lie(Spin(7)) = \set{\lambda^2(X) | X \in \lieso(7)}.
\end{equation}

Let $\set{e_i \pm e_j | 1\leq i < j \leq 4}$ be the positive root system of $\lieso(8)$ with
vanishing $e_4$ on $\lieso(7)$, then the automorphism $\lambda^2$ induces the following
transformation of $e_i$.
\begin{equation}\label{matrixlambda}
\lambda^2
\begin{pmatrix}
e_1\\
e_2\\
e_3\\
e_4
\end{pmatrix}
=\left(
\begin{array}{rrrr}
\half & \half & \half & \half \\
\half & \half & -\half & -\half \\
\half & -\half & \half & -\half \\
-\half & \half & \half & -\half \\
\end{array}\right)
\begin{pmatrix}
e_1\\
e_2\\
e_3\\
e_4
\end{pmatrix}.
\end{equation}

\smallskip

In (\ref{eqnspin7}) the $\theta_2$-component acts transitively on $\sph^7 \subset \Cay$ and the
isotropy subgroup at $v_0$ is $G_2$. Since $\kappa$ is the identity map on $SO(7)$,
$\theta_1(x)\theta_2(y) = \theta_2(xy)$ for $\theta_2\in G_2\subset SO(7)$ which implies $\theta_1
= \theta_2$. Hence we have
\[
G_2 = \set{(\theta, \theta)| \theta \in SO(7) \mbox{ and }\theta(x)\theta(y) = \theta(xy), \mbox{
for any } x, y \in \Cay},
\]
and the Lie algebra is
\[
\frakg_2 = \set{(X,X)| X\in \lieso(7), X = \lambda(X)}.
\]

If we identify $\Lie(Spin(7))$ with $\lieso(7)$, then we have
\begin{equation*}
\frakg_2 = \set{X | X \in \lieso(7), X=\lambda(X)}.
\end{equation*}
Choose a basis of $\Cay$ over $\Real$, then we can write down the explicit embedding of $\frakg_2$,
see also in \cite{Murakami}. The typical element of $\frakg_2$ has the following form:
\begin{equation*}
X=
\begin{pmatrix}
0 & x_1 - y_1 & x_3 + y_3 & -x_2 + y_2 & -x_4 - y_4 & x_6 + y_6 & x_5 -y_5 \\
-x_1 + y_1 & 0 & a & y_5 & y_6 & y_4 & y_2 \\
-x_3 - y_3 & -a & 0 & x_6 & x_5 & x_2 & x_4 \\
x_2 -y_2 & - y_5 & -x_6 & 0 & b & y_3 & y_1 \\
x_4+y_4 & - y_6 & -x_5 & -b & 0 & x_1 & x_3 \\
-x_6 - y_6 & - y_4 & -x_2 & - y_3 & -x_1 & 0 & a+b \\
-x_5+y_5 & - y_2 & -x_4 & - y_1 & -x_3 & -a-b & 0
\end{pmatrix},
\end{equation*}
where $a, b, x_1, \cdots, x_6, y_1, \cdots, y_6 \in \Real$.

The linear functionals $\set{e_1, e_2, e_3}$ of the Cartan subalgebra of $\lieso(7)$ satisfy the
relation $e_1= e_2 + e_3$ when restricted to the Cartan subalgebra of $\frakg_2$. Suppose
$\set{\alpha_1, \alpha_2}$ be the set of the positive simple roots of $G_2$ where $\alpha_1$ is the
short one, then under the restriction, $\overline{e_3} = \alpha_1$ and $\overline{e_2 - e_3} =
\alpha_2$.

\smallskip

We consider the last pair $(G_2, SU(3))$. $G_2$ acts transitively on the unit sphere $\sph^6 =
\set{x \in \Cay | \norm{x} = 1, <x, v_0> = 0}$. Let $v_1$ be a unit element which is orthogonal to
$v_0$, then the isotropy group at $v_1$ is isomorphic to $SU(3) \subset SO(6)$. Since $G_2$ and
$SU(3)$ share the same maximal torus, the restriction of the roots of $\frakg_2$ is the identity
map.

\smallskip

In the proof we will use the classical branching rules for the special orthogonal groups which are
well-known and the proof can be found, for examples, in \cite{Knappbook}.

\begin{thm}[Branching Rule for $(\lieso(2k+1), \lieso(2k))$]\label{thmBranchingsoodd}
The irreducible representation with highest weight $a_1e_1 + a_2e_2 + \cdots + a_ke_k$ of
$\lieso(2k+1)$ decomposes with multiplicity $1$ under $\lieso(2k)$, and the representations of
$\lieso(2k)$ that appear are exactly those with highest weights $\mu = c_1e_1+ c_2e_2+ \cdots
+c_ke_k)$ such that
\[
a_1 \geq c_1 \geq a_2 \geq c_2 \geq \cdots \geq a_k \geq \abs{c_k},
\]
where $a_i (i = 1,2,\cdots,k)$ are integers or all half integers and $c_j (j = 1,2,\cdots,k)$ are
all integers or all half integers.
\end{thm}

\begin{thm}[Branching Rule for $(\lieso(2k), \lieso(2k-1))$]\label{thmBranchingsoeven}
The irreducible representation with highest weight $a_1e_1 + a_2e_2 + \cdots + a_ke_k$ of
$\lieso(2k)$ decomposes with multiplicity $1$ under $\lieso(2k-1)$, and the representations of
$\lieso(2k-1)$ that appear are exactly those with highest weights $\mu =c_1e_1+ c_2e_2+ \cdots
+c_{k-1}e_{k-1})$ such that
\[
a_1 \geq c_1 \geq a_2 \geq c_2 \geq \cdots \geq c_{k-1} \geq \abs{a_k},
\]
where $a_i (i = 1,2,\cdots,k)$ are integers or all half integers and $c_j (j = 1,2,\cdots,k-1)$ are
all integers or all half integers.
\end{thm}

\smallskip

\emph{Proof of the 3rd Part of Theorem \ref{thmclassonecpx}:} First let $K = Spin(9)$, $H=Spin(7)$
and $\frakk = \Lie(Spin(9))$, $\frakh = \Lie(Spin(7))$ be their Lie algebras. $L=Spin(8)$ lies
between them and denote its Lie algebra $\Lie(Spin(8))$ by $\frakl$. Let
\[
\Delta_K^+ = \set{e_i \pm e_j| 1\leq i < j \leq 4} \cup \set{e_i | 1\leq i \leq 4},
\]
be the positive root system of $\frakk$, and then $\frakl$ has the positive roots system
\[
\Delta_L^+=\set{e_i \pm e_j| 1\leq i < j \leq 4}.
\]
Since $\frakk$ and $\frakl$ share the same Cartan subalgebra, $\Sigma=\set{e_1, e_2, e_3, e_4}$.
Therefore the branching rule for the pair $(\frakk, \frakl)$ is the same as the classical one for
$(\lieso(9), \lieso(8))$.

Let $f_i = \lambda^2(e_i)$ for $i=1,2,3,4$. Since $\lambda$ is an automorphism of $\lieso(8)$, we
can write
\[
\Delta^+_L =\set{f_i \pm f_j| 1\leq i < j \leq 4},
\]
and $\frakh$ has the positive root system
\[
\Delta^+_H = \set{f_i \pm f_j| 1\leq i < j \leq 3} \cup \set{f_i | 1\leq i \leq 3}.
\]
Therefore $\Sigma=\set{f_1, f_2, f_3}$. The Weyl group $W_K$ also acts on $f_i$'s and $\delta_K =
3e_1 + 2e_2 + e_3 = 3f_1 + 2f_2 + f_3$. Suppose $\mu$ be an irreducible representation of $K$ with
the highest weight $\mu = b_1f_1+ b_2f_2+b_3f_3+b_4f_4$, then by the classical branching rule for
$(\lieso(8), \lieso(7))$, $\mu$ decomposes with multiplicity $1$ when restricted to $\frakh$, and
the representations of $\frakh$ that appear are exactly those with highest weights $\tau =c_1f_1+
c_2f_2+ c_3f_3$ such that
\[
\begin{array}{l}
b_1 \geq c_1 \geq b_2 \geq c_2 \geq  b_3 \geq c_3 \geq \abs{b_4},
\end{array}
\]
where $b_i (i = 1,2,3,4)$ and $c_j (j = 1,2,3)$ are all integers or half integers. Therefore the
class one representations of the pair $(K,H)$ are exactly those with the highest weights $bf_1$
where $b \in \Int$, or $\half b(e_1 + \cdots + e_4)$ and the multiplicity of trivial representation
is $1$.

Now suppose that $\mu$ is a class one representation of the pair $(K,H)$ with the highest weight
$\rho = a_1 e_1 + \cdots + a_4 e_4$, then $\Res^K_L \rho$ contains at least one representation with
highest weight $\half b (e_1 + \cdots + e_4)$. From the branching rule of $(\lieso(9), \lieso(8))$,
we have
\[
a_1 \geq \half b \geq a_2 \geq \half b \geq a_3 \geq \half b \geq a_4 \geq \half b.
\]
Therefore $a_1 \geq a_2 = a_3 = a_4 \geq 0$ and all are integers or half integers. If $a_i$'s
satisfy the conditions, then $\rho$ contains only one representation with the highest weight $\mu$
which implies that $[\Res^K_H \rho : Id] =1$.

Next we look at the pair $(K,H) = (Spin(7), G_2)$. We identify $\frakk$ with $\lieso(7)$ so that
\[
\Delta_K^+ = \set{e_i \pm e_j| 1\leq i < j \leq 3} \cup \set{e_i | 1\leq i \leq 3}.
\]
By the relation $e_1 = e_2 + e_3$ and the restriction $\overline{e_3} = \alpha_1$, $\overline{e_2 -
e_3} = \alpha_2$, we have
\[
\overline{\Delta_K^+}=\set{\alpha_1, 3\alpha_1+2\alpha_2, \alpha_1 + \alpha_2, 3\alpha_1 +
\alpha_2,\alpha_2, 2\alpha_1+\alpha_2, 2\alpha_1+\alpha_2, \alpha_1+\alpha_2, \alpha_1}
\]
under the restriction. In terms of $\alpha_1, \alpha_2$, $\frakh$ has the following positive root
system
\[
\Delta_H^+ = \set{\alpha_1, \alpha_2, \alpha_1+ \alpha_2, 2\alpha_1 + \alpha_2, 3\alpha_1 +
\alpha_2, 3\alpha_1 + 2\alpha_2}.
\]
Therefore $\Sigma = \set{\alpha_1, \alpha_1 + \alpha_2, 2\alpha_1 + \alpha_2}$. The Weyl group
$W_K$ permutes $\set{e_1, e_2, e_3}$ and changes the sign of each $e_i$.

Let $\rho$ be a class one representation with the highest weight $\rho = a_1e_1 + a_2e_2 + a_3e_3$
$(a_1 \geq a_2 \geq a_3 \geq 0)$. Then from the Kostant's Branching Theorem \ref{thmKostant}, by a
computation, the multiplicity of the trivial representation $\Id$ of $G_2$ in $\Res \rho$ is
\begin{equation*}
p = (a_1 + a_3 + 1) - (a_1 -a_3) - (a_1 + a_3) + \max\{0, a_1 - a_3 -1\}-(a_1 - a_3 -1) + \max\{0,
a_1 -a_3 -1\}.
\end{equation*}
Therefore $p\neq 0$ if and only if $a_1 = a_3$, i.e., $a_1 = a_2 = a_3$. Thus in this case we have
$p=1$.

Finally we look at the pair $(K, H) = (G_2, SU(3))$.  $\frakk=\frakg_2$ and $\frakh=\liesu(3)$
share the same Cartan subalgebra and hence the restriction of the roots is the identity map.
$\frakk$ has the following positive root system
\[
\Delta_K^+ = \set{\alpha_1, 3\alpha_1 + \alpha_2, 2\alpha_1 + \alpha_2, 3\alpha_1 + 2\alpha_2,
\alpha_1 + \alpha_2, \alpha_2},
\]
where $\alpha_1= e_3$, $\alpha_2 = e_2 - e_3$ and $e_2 + e_3 -e_1 = 0.$ Among the $12$ roots of
$\frakk$, the six long roots consist the root system of $\frakh$, i.e.
\[
\Delta_H^+ = \set{3\alpha_1 + 2\alpha_2, 3\alpha_1 + \alpha_2, \alpha_2}.
\]
Therefore $\Sigma= \set{\alpha_1, \alpha_1 + \alpha_2, 2\alpha_1 + \alpha_2}.$ The Weyl group $W_K$
is the symmetry group of the regular hexagon $D_6 = <\sigma, \tau | \sigma^6 = 1, \tau^2 = 1>$ with
the following operations on $\alpha_1, \alpha_2$:
\[
\begin{array}{rcl}
\sigma & : & \alpha_1 \mapsto -\alpha_1 - \alpha_2, \quad \alpha_2 \mapsto 3\alpha_1 + 2\alpha_2,\\
\tau & : & \alpha_1 \mapsto -\alpha_1 -\alpha_2, \quad \alpha_2 \mapsto \alpha_2.
\end{array}
\]

Suppose $\rho$ is a class one representation with the highest weight $\rho = a_1\fw_1 + a_2\fw_2 =
(2a_1 + 3a_2)\alpha_1 + (a_1+2a_2)\alpha_2$, where $a_1,a_2$ are two nonnegative integers. A
computation using Kostant Branching Theorem shows that $a_2 =0$ and the multiplicity $[\Res \rho :
\Id]$ is $1$. \qedbox

\smallskip

\smallskip

\subsection{Properties of Class One Representations}\label{secproperties}

First we consider the types and nontrivial kernels of class one representations. They follow easily
from Proposition 23.13 and propositions in $\S$26.3 in \cite{FultonHarris}.

\begin{prop}\label{propclassonetypekernel}
The type and non-trivial kernel $Z$ of each complex class one representation is listed in Table
\ref{tabletypeclassone}.

\begin{table}[!h]
\begin{center}
\begin{tabular}{|c|c|l|l|}
  \hline
  $K$  & $\rho$  &  $\quad \quad \quad$ Type  & $\quad\quad\quad\quad\quad$ $Z$ \\
  \hline
  \hline
  $SO(n)$ & $a \fw_1$ & real  &  $\Int_2$: if both $n$ and $a$ are even \\
  \hline
  $SU(n)$ & $a\fw_1 + b \fw_{n-1}$ & real: if $a=b$ & $\Int_l$, $l = \gcd (a+b, n)$ \\
          &                        & complex: otherwise & \\
  \hline
  $U(n) $ & $ae_1 -b e_n +$ & real: if $a=b$ &$\Int_{\abs{(a-b)(1+mn)}}$: if $a\ne b$ \\
          & $m(a-b)(e_1 + \cdots + e_n)$ & complex: otherwise & \\
  \hline
  $Sp(n)$ & $a\fw_1 + b\fw_2$ & real: $a$ is even & $\Int_2$: if $a$ is even \\
          &                       & quaternionic: otherwise &\\
  \hline
  \hline
  $Sp(n)\times Sp(1)$ & $(a \fw_1 + b \fw_2)\otimes a \fw_1$ & real & $\Int_2 \times \Int_2$: if $a\ne 0$ is even \\
                      &                                      &      & $\Int_2 \times Sp(1)$: if
                      $a=0$ \\
  \hline
                     &                                     & complex: $k\ne 0$ & $\Id \times \Int_{|k|}$: if $a$ is
                     odd \\
  $Sp(n)\times U(1)$ & $(a \fw_1 + b\fw_2) \otimes \phi^k$ & real: $k=0$ and $a$ is even & $\Int_2 \times \Int_{|k|}$: if $a$ is
                                                                                  even\\
                     &                                     & quaternionic: otherwise &\\
  \hline
  \hline
  $G_2$ & $a \fw_1$ & real & $\quad \quad$---------------------  \\
  \hline
  $Spin(7) $ & $a\fw_3$ & real & $\Int_2$: if $a$ is even \\
  \hline
  $Spin(9) $ & $a\fw_1 + b\fw_4$ & real & $\Int_2$: if $b$ is even \\
  \hline
  \hline
  $U(1)$ & $\phi^k$ & complex & $\Int_{\abs{k}}$: if $\abs{k}>1$ \\
  \hline
\end{tabular}
\vskip 0.1cm \caption{Type and non-trivial kernel of class one representations}
\label{tabletypeclassone}
\end{center}
\end{table}
\end{prop}

In Table \ref{tabletypeclassone}, for $K=U(n)$($a=b$) and $K=Sp(n)\times U(1)$($k=0$), the $\Int_0$
factor in the non-trivial kernel should be interpreted as the group $U(1)=\set{z\cdot \Id | z \in
U(1)}$. Furthermore, $\Int_2 \subset K$ is generated by $-\Id$ and $\Int_l \subset SU(n)$(or
$U(n)$) is generated by $\varphi_l \cdot \Id$ with $\varphi_l = \exp(\frac{2\pi\sqrt{-1}}{l})$. The
notation $\gcd(p,q)$ stands for the greatest common divisor of $p$, $q$.

\smallskip

Next we give the dimensions of the class one representations, see \cite{FultonHarris}.
\begin{prop}\label{propclassonedim}
The dimension of each complex class one representation for the pair $(K,H)$ is listed as follows:
\begin{itemize}
  \item $(SO(n),SO(n-1))$ :  $\rho = a \fw_1$, \\
  $\dim \rho = \frac{2a+n-2}{a+n-2}\binom{a+n-2}{a}$ \\
  \item $(SU(n),SU(n-1))$ :  $\rho = a \fw_1 +b \fw_{n-1}$, \\
  $\dim \rho = \frac{a+b+n-1}{n-1}\binom{a+n-2}{a} \binom{b+n-2}{b}$ \\
  \item $(U(n),U(n-1)_m)$ :  $\rho = a e_1 - b e_n + m(a-b)(e_1 + \cdots + e_n)$, \\
  $\dim \rho = \frac{a+b+n-1}{n-1}\binom{a+n-2}{a} \binom{b+n-2}{b}$ \\
  \item $(Sp(n),Sp(n-1))$ :  $\rho = a \fw_1 + b\fw_2$, \\
  $\dim \rho = \frac{(a+1)(a+2b+2n-1)}{(a+b+1)(a+b+2n-1)}  \binom{a+b+2n-1}{a+b} \binom{b+2n-3}{b}$
  \\
  \item $(Sp(n)\times Sp(1),Sp(n-1)\times Sp(1))$ :  $\rho = (a \fw_1 + b\fw_2)\otimes a\fw_1$, \\
  $\dim \rho = \frac{(a+1)^2 (a+2b+2n-1)}{(a+b+1)(a+b+2n-1)}  \binom{a+b+2n-1}{a+b} \binom{b+2n-3}{b}$
  \\
  \item $(Sp(n)\times U(1) ,Sp(n-1)\times U(1)_m)$ :  $\rho = (a \fw_1 + b\fw_2)\otimes \phi^k $, \\
  $\dim \rho = \frac{(a+1)(a+2b+2n-1)}{(a+b+1)(a+b+2n-1)}  \binom{a+b+2n-1}{a+b} \binom{b+2n-3}{b}$
  \\
  \item $(G_2,SU(3)) $    :  $\rho = a\fw_1$, \\
  $\dim \rho = \frac{1}{120}(a+1)(a+2)(a+3)(a+4)(2a+5)$ \\
  \item $(Spin(7),G_2)$   :  $\rho = a\fw_3 $, \\
  $\dim \rho = \frac{1}{360}(a+1)(a+2)(a+3)^2(a+4)(a+5)$ \\
  \item $(Spin(9),Spin(7))$: $\rho = a\fw_1 + b\fw_4$, \\
  $\dim \rho = \frac{1}{1814400}(a+1)(a+2)(a+3)(b+1)(b+2)(b+3)^2(b+4)(b+5)(a+b+4)(a+b+5)(a+b+6)(2a+b+7)
  $ \\
  \item $(U(1),\set{1})$: $\rho = \phi^k$, $\dim \rho =1$.
\end{itemize}
\end{prop}

\smallskip

Finally for each spherical pair $(K,H)$, we consider the embedding $\rho: K \To SO(l)$, $U(l)$ and
$Sp(l)$. We compare the dimension of the sphere $K/H$ and $l$.

\begin{prop}\label{propdimcomparereal}
Table \ref{tableclassonedimcomparereal} is the list of the dimension $s$ of the sphere $K/H$ and
the embedding $\rho: K \To SO(l)$ with $l\leq k_0(s+1)$ where $k_0$ is equal to $1$, $2$ or $4$ if
$\rho$ is of real, complex or quaternionic type.
\begin{table}[!h]
\begin{center}
\begin{tabular}{|c|c|c|c|c|c|}
  \hline
  $K$  & $H$ & $\rho$  &  $l$  & $s+1$ & $k_0(s+1)$ \\
  \hline
  \hline
  $SO(n)$ & $SO(n-1)$ & $\fw_1$ & $n$ & $n$ & \\
  \hline
  $SU(n)$ & $SU(n-1)$ & $[\fw_1]_\Real$ & $2n$ & $2n$ & $4n$ \\
  \hline
  $U(n) $ & $U(n-1)_m$  & $[e_1 + m(e_1 + \cdots + e_n)]_\Real$ & $2n$ & $2n$ & $4n$ \\
  \hline
  $Sp(n)$ & $Sp(n-1)$ & $[\fw_1]_\Real$ & $4n$ & $4n$ & $16n$\\
  \hline
  \hline
  $Sp(n)\times Sp(1)$ & $Sp(n-1)\times Sp(1)$ & $\fw_1 \otimes  \fw_1$ & $4n$ & $4n$ &  \\
  \hline
  $Sp(n)\times U(1)$ & $Sp(n-1) \times U(1)_m$ & $[\fw_1 \otimes  \phi^k]_\Real$ & $4n$ & $4n$ & $8n$ or $16n$ \\
  \hline
  \hline
  $G_2$ & $SU(3)$ & $\fw_1$ & $7$ & $7$ & \\
  \hline
  $Spin(7) $ & $G_2$ & $\fw_3$ & $8$ & $8$ & \\
  \hline
  $Spin(9) $ & $Spin(7)$ & $\fw_1$, $\fw_4$ & $9$, $16$ & $16$ &  \\
  \hline
  \hline
  $U(2)$ & $U(1)_m$ & $e_1 - e_2$ & $3$ & $4$ &\\
  \hline
  $Sp(2)   $ & $Sp(1)$ & $\fw_2$ & $5$ & $8$ & \\
  \hline
  $Sp(1)$ & $\set{1}$ & $2\fw_1$ & $3$ & $4$ & \\
  \hline
  \hline
  $U(1)$ & $\set{1}$ & $[\phi^k]_\Real$ & $2$ & $2$ & $4$ \\
  \hline
  $U(2)$ & $U(1)_m$ & $[2e_1 + 2m(e_1+e_2)]_\Real$ & $6$ & $4$ & $8$ \\
         &          & $[3e_1 + 3m(e_1+e_2)]_\Real$ & $8$ & &  \\
  \hline
  $U(3)$ & $U(2)_m$ & $[2e_1 + 2m(e_1 + e_2 + e_3)]_\Real$ & $12$ & $6$ & $12$ \\
  \hline
  $SU(3)$ & $SU(2)$ & $[2\fw_1]_\Real$ & $12$ & $6$ & $12$ \\
  \hline
  \hline
  $Sp(2)$ & $Sp(1)$ & $[\fw_1 + \fw_2]_\Real$ & $32$ & $8$ & $32$ \\
  \hline
  $Sp(1)$ & $\set{1}$ & $[3\fw_1]_\Real$, $[5\fw_1]_\Real$, $[7\fw_1]_\Real$ & $8,12,16$ & $4$ & $16$ \\
  \hline
\end{tabular}
\vskip 0.1cm \caption{Orthogonal Class One Representation with Small Dimension}
\label{tableclassonedimcomparereal}
\end{center}
\end{table}
\end{prop}

For the complex representations, we have

\begin{prop}\label{propdimcomparecpx}
Table \ref{tableclassonedimcomparecpx} is the list of the dimension $s$ of the sphere $K/H$ and the
embedding $\rho: K \To U(l)$ with $2l\leq s+1$.
\begin{table}[!h]
\begin{center}
\begin{tabular}{|c|c|c|c|c|}
  \hline
  $K$  & $H$ & $\rho$  &  $l$  & $s+1$ \\
  \hline
  \hline
  $SU(n)$ & $SU(n-1)$ & $\fw_1$, $\fw_{n-1}$ & $n$, $n$ & $2n$ \\
  \hline
  $U(n) $ & $U(n-1)$  & $e_1$, $-e_n$ & $n$, $n$ & $2n$ \\
  \hline
  $Sp(n)$ & $Sp(n-1)$ & $\fw_1$ & $2n$ & $4n$\\
  \hline
  $Sp(n)\times U(1)$ & $Sp(n-1) \times U(1)_m$ & $\fw_1 \otimes  \phi^k$ & $2n$ & $4n$ \\
  \hline
\end{tabular}
\vskip 0.1cm \caption{Complex Class One Representation with Small Dimension}
\label{tableclassonedimcomparecpx}
\end{center}
\end{table}
\end{prop}

A similar result for the quaternionic representations is

\begin{prop}\label{propdimcomparequa}
The quaternionic class one representations $\rho: K \To Sp(l)$ with $4 l \leq \dim (K/H) + 1$ are
exactly the representation $\rho = \fw_1$ for the pair $(Sp(n),Sp(n-1))$($n\geq 1$).
\end{prop}

\medskip

\medskip


\vfill
\end{document}